\documentclass[12pt,leqno]{article}
\usepackage{amsmath}
\usepackage{amssymb}
\usepackage{amscd}
\renewcommand{\cite}[1]{[#1]}
\def\overset#1#2{{\mathrel{\mathop {{#2}_{}}\limits^{#1}}}}
\def\underset#1#2{{\mathrel{\mathop {{}_{} {#2}}\limits_{{#1}_{}}}}}
\def\upplim_#1{\underset{#1}{\overline\lim}\;}
\def\lowlim_#1{\underset{#1}{\underline\lim}\;}
\newtheorem{claim}[equation]{\indent \rm Claim}%

\newtheorem{corollary}[equation]{Corollary}
\newtheorem{definition}[equation]{\indent{\it Definition}\rm }

\newtheorem{lemma}[equation]{Lemma}
\newtheorem{lem}[equation]{Lemma}
\newtheorem{proposition}[equation]{Proposition}

\newtheorem{convention}[equation]{\indent \rm {\it Convention}}
\newtheorem{theorem}[equation]{Theorem}

\newcommand{\C}{{\mathbf{C}}}
\newcommand{\codim}{{\mathrm{codim}}}

\newcommand{\del}{{\partial}}
\newcommand{\delbar}{\bar{\partial}}
\newcommand{\I}{\mathcal{I}}

\newcommand{\Lie}{{\mathop{\mathrm{Lie}}}}
\newcommand{\N}{\mathbf{N}}
\newcommand{\NO}{NO$\frac{84}{90}$}
\renewcommand{\O}{{\mathcal{O}}}
\newcommand{\ord}{{\mathrm{ord}}}
\renewcommand{\P}{{\mathbf{P}}}
\newcommand{\pnc}{{\mathbf{P}^n(\mathbf{C})}}
\newcommand{\pNc}{{\mathbf{P}^N(\mathbf{C})}}
\newcommand{\Pic}{{\mathop{\mathrm{Pic}}}}

\newcommand{\R}{{\mathbf{R}}}

\newcommand{\Sing}{\mathrm{Sing}}
\newcommand{\St}{\mathrm{St}}
\newcommand{\StD}{\mathrm{St}(D)}
\newcommand{\supp}{\mathrm{Supp}\,}
\newcommand{\tensor}{\otimes}

\newenvironment{proof}{\par{\it Proof.}}{{\it Q.E.D.}\par\vskip3pt}
\numberwithin{equation}{section}

\title{
The Second Main Theorem for Holomorphic\\
Curves into Semi-Abelian Varieties II
\thanks{  Research supported in part by Grant-in-Aid
   for Scientific Research (A)(1) 13304009
and (S) 17104001.}}
\author{
Junjiro Noguchi, J\"org Winkelmann and Katsutoshi Yamanoi}
\begin{document}
\setlength{\baselineskip}{18pt}
\maketitle
\begin{abstract}
We establish the second main theorem with the best truncation level
one
$$
T(r; \omega_{\bar Z, J_k(f)})\leqq N_1(r; J_k(f)^* Z) +\epsilon T_f(r)||_\epsilon
$$
for the $k$-jet lift $J_k(f):\C \to J_k(A)$ of
an algebraically non-degenerate entire holomorphic curve
$f:\C \to A$ into a semi-abelian variety $A$
and an arbitrary algebraic reduced subvariety $Z$ of $J_k(A)$;
the low truncation level is important for applications.
Finally we give some applications, including the solution of a problem
posed by Mark Green (1974).
\end{abstract}

\section{Introduction and main result}

Let $f: \C \to V$ be a holomorphic curve into a complex projective
manifold $V$ with Zariski dense image
and let $D$ be an effective reduced divisor on $V$.
Under some ampleness condition for the space
$H^0(V, \Omega^1_V(\log D))$ of logarithmic 1-forms along $D$
we proved in \cite{N77}, \cite{N81} the following inequalities of
the second main theorem type,
\begin{align*}
\kappa T_f(r) &\leqq N (r; f^*D) + O(\log r) + O(\log T_f(r)) ||, \\
\kappa' T_f(r) &\leqq N_1 (r; f^*D) + O(\log r) + O(\log T_f(r)) ||,
\end{align*}
where $T_f(r)$ denotes the order function of $f$,
$N(r; f^*D)$ (resp.\ $N_l(r; f^*D)$) the counting function
(resp.\ truncated to level $l$) of the pull-backed divisor $f^*D$,
and $\kappa$ and $\kappa'$ are positive constants (cf.\ \S2).
It is an interesting and fundamental problem to determine
the constant $\kappa$ or $\kappa'$.
In the case where $V$ is the compactification of a semi-abelian
variety $A$ this
problem is related to what kind of compactification $V$ of
$A$ we take.
In our former paper \cite{NWY02} we proved that for a holomorphic
curve $f:\C \to A$ into a semi-abelian variety $A$
and an algebraic divisor $D$ on $A$,
\begin{equation}
\label{0.1}
T_f(r; L(\bar D)) \leqq N_l(r; f^*D)+ O(\log r) +
O(\log T_f(r; L(\bar D)))||.
\end{equation}
Here we  used a compactification $\bar A$ of $A$ such that
the maximal 
affine subgroup $(\C^*)^t$ of $A$ was compactified by $(\P^1(\C))^t$,
and we assumed a boundary condition (Condition 4.11
in \cite{NWY02}) for the closure $\bar D$ of $D$
in $\bar A$;
this roughly meant the divisor
$\bar D + (\bar A \setminus A)$ to be in general position
and has been expected to be removed by a suitable choice of
a compactification of $A$.
It is an important and very interesting problem to take
the truncation level
$l$  as small as possible.

Let $X_k(f)$ denote the Zariski closure of the image of the $k$-jet
lift of $f$ in the $k$-jet space $J_k(A)$ over $A$.
The purpose of this paper is to prove (cf.\ \S\S 2, 3 for notation)
\medskip

{\bf Main Theorem.}  {\it Let $A$ be a semi-abelian variety.
Let $f: \C \to A$ be a holomorphic curve with Zariski
dense image.

{\rm (i)}  Let $Z$ be an algebraic reduced subvariety of $X_k(f)$
$(k \geqq 0)$.
Then there exists a compactification $\bar X_k(f)$ of $X_k(f)$
such that
\begin{equation}
\label{1.1}
T(r; \omega_{\bar Z , J_k(f)}) \leqq N_1 (r; J_k(f)^* Z)+
\epsilon T_f(r) ||_\epsilon,\quad \forall \epsilon >0,
\end{equation}
where $\bar Z$ is the closure of $Z$ in $\bar X_k(f)$.

{\rm (ii)}  Moreover, if $\codim_{X_k(f)}Z \geqq 2$, then
\begin{equation}
\label{1.2}
T(r; \omega_{\bar Z, J_k(f)}) \leqq
\epsilon T_f(r) ||_\epsilon,\quad \forall \epsilon >0.
\end{equation}

{\rm (iii)}  In the case when $k=0$ and $Z$ is an effective divisor
$D$ on $A$, the compactification $\bar A$ of $A$
can be chosen as smooth, equivariant with respect to the $A$-action,
and independent of $f$; furthermore, \eqref{1.1} takes the form
\begin{equation}
\label{1.3}
T_f(r; L(\bar D)) \leqq N_1(r; f^*D)
+\epsilon T_f(r; L(\bar D))||_\epsilon,
\quad \forall \epsilon>0.
\end{equation}
}

Note that in the above estimate \eqref{1.1},
\eqref{1.2} or \eqref{1.3}
 the small error term ``$\epsilon T_f(r)$''
cannot be replaced by ``$O(\log r) + O(\log T_f(r))$''
(see \cite{NWY02} Example (5.36)).

The Main Theorem is an advancement of \cite{NWY02} and \cite{Y04}.
When $A$ is an abelian variety,
\eqref{1.3}  was proved by Yamanoi \cite{Y04}
(cf.\ \cite{Y04} (3.1.8)).
A key of the proof of \eqref{0.1} in \cite{NWY02}
was Lemma 5.6 at p.~147, and here we will again use the same idea
for jets of jets (see Claim \ref{jetdiff3}).

There is a related result due to Siu-Yeung \cite{SY03},
where they obtained \eqref{0.1} with an improved truncation level $l=l(D)$
dependent only on the Chern numbers of $D$.
In their proof the key was Claim 1
at p.~443 which was the same as \cite{NYW02} Lemma 5.6 restricted to
the abelian case with a computation of intersection numbers.

It is interesting to observe that the error term being
``$ O(\log r) +O(\log T_f(r; L(\bar D)))||$'',
the truncation level $l$ in \eqref{0.1} has to depend on $D$,
but the error term being allowed to be a littel bit large,
``$\epsilon T_f(r; L(\bar D)))||_\epsilon$'',
$l$ can be one, the smallest possible.
In applications, the truncation of level one is very definite.

To deal with semi-abelian varieties the main difficulties are
caused by the following two points:
\begin{enumerate}
\item
Semi-abelian varieties are not compact and need some good
     compactifications.
\item
There is no Poincar\'e reducibility theorem for semi-abelian varieties.
\end{enumerate}
It is also noted that a part of the proof of the Main Theorem
for abelian varieties in \cite{Y04} does not hold for semi-abelian
varieties (\cite{Y04} \S3 Claim), and that a different
and considerably simpler proof for that part will be provided
(see Lemma \ref{codim2}).

In \S7 we will give two applications of the Main Theorem.
The first is a complete affirmative answer to a conjecture
of M. Green \cite{G74} pp. 229--230
(cf.\ Theorem \ref{greenconj}).
The second is a non-existence theorem for some differential equations
defined over semi-abelian varieties (cf.\ Theorem \ref{deq}).

More applications will be obtained in \cite{NWY05}.

{\it Acknowledgement.}   We learned the conjecture of
M. Green \cite{G74} from Professor A.E. Eremenko, to whom
we are very grateful.

\section{Notation}

The notation here follows that of [NWY02].
For a general reference of this section, cf.\ [\NO].
For convenience we recall some of definitions.
Let $M$ be a compact complex manifold and let
$\omega$ be a smooth (1,1)-form on $M$.
Let $f:\C \to M$ be a holomorphic curve into $M$.
We define the order function of $f$ with respect to $\omega$ by
\begin{equation}
\label{2.1}
T_f(r; \omega)=\int_1^r\frac{dt}{t}\int_{|z|<t} f^*\omega
\qquad (r>1).
\end{equation}
If $M$ is K\"ahler and $d\omega=0$,
$$
T_f(r; \omega)=T_f(r; \omega')+O(1)
$$
for a $d$-closed (1,1)-form $\omega'$ in the same cohomology class
$[\omega] \in H^2(M, \R)$.
Therefore we set, up to $O(1)$-term,
\begin{equation}
\label{2.2}
T_f(r; [\omega])=T_f(r; \omega).
\end{equation}
Let $L \to M$ be a hermitian line bundle with Chern class $c_1(L)$.
Then we set
$$
T_f(r; L)=T_f(r; c_1(L)),
$$
which is defined again up to $O(1)$-term.

For a divisor $D$ on $M$ we denote by $L(D)$ the line bundle
determined by $D$.

Let $E=\sum_{\mu=1}^\infty \nu_\mu z_\mu$ be a divisor on $\C$
with distinct $z_\mu \in \C$.
Then we set
$$
\ord_z E =\begin{cases}
\nu_\mu, & z=z_\mu, \\
0, & z \not\in \{z_\mu\}.
\end{cases}
$$
We define the counting functions of $E$ truncated to $l \leqq \infty$ by
\begin{align*}
n_l(t; E) &=\sum_{\{|z_\mu|<t\}} \min \{\nu_\mu, l\},\\
N_l(r; E) &=\int_1^r \frac{n_l(t; E)}{t}dt.
\end{align*}
We define the counting functions of $E$ by
$$
n(t;E)=n_\infty(t; E), \qquad N(r; E)=N_\infty(r; E).
$$

{\it Definition of small terms.}
(i)  For a line bundle $L \to M$ and
a holomorphic curve
$f:\C \to M$ we denote by $S_f(r; L)$
such a small term as
$$
S_f(r; L)= O(\log r) + O(\log^+ T_f(r; L))||,
$$
where ``$||$'' stands for the inequality to hold
for every  $r>1$ outside a Borel set of finite Lebesgue measure.

(ii)  Let $h(r)$ $(r>1)$ be a real valued function. We write
$$
h(r) \leqq \epsilon T_f(r; L)||_\epsilon,\quad \forall \epsilon>0,
$$
if the stated inequality holds for every  $r>1$ outside
a Borel set of finite Lebesgue measure, dependent on an
arbitrarily given $\epsilon >0$.

{\it Definition.}  When $M$ is an algebraic variety,
we say that $f:\C \to M$
is {\it algebraically (resp.\ non-) degenerate} if the image
$f(\C)$ is (resp.\ not) contained in a proper algebraic subset
of $M$.

The following follows from general properties of
order functions (\cite{\NO}).

\begin{lemma}
\label{2.4}
Let $f:\C \to M$ be a holomorphic curve into a complex projective
manifold $M$ and  $H$ a line bundle on $M$.
Assume that $H$ is big,
and that $f$ is algebraically non-degenerate.
Then 
\[
T_f(r,L)= O(T_f(r,H))
\]
for every line bundle $L$ on $M$.
\end{lemma}

If $f:\C\to M$ is algebraically degenerate, we may consider
the Zariski closure $N$ of $f(\C)$ and a desingularization
$\tau:\tilde N\to N$. Then $f$ lifts to a map to $\tilde N$ and
$\tau^*(H|_N)$ is big on $\tilde N$ for every ample line bundle
$H$ on $M$.
As a consequence we obtain:

\begin{lemma}
\label{2.4a}
Let $f:\C \to M$ be a holomorphic curve into a complex projective
manifold $M$.
Let $h(r)$ be a non-negative valued function in $r>1$.
Then $h(r)=S_f(r; H)$ holds for every ample line bundle
if and only if it holds for at least one ample line bundle.

Similarily the statement
$h(r) \leqq \epsilon T_f(r; H)||_\epsilon, \forall \epsilon>0$ ,
respectively
$h(r)=O(T_f(r; H))$ 
holds for every ample line bundle $H$ if and only if 
it holds for at least one ample line bundle.
\end{lemma}

If one of these conditions
holds for one and therefore for all ample line bundles $H$,
we simply write $h(r)=S_f(r)$ (resp.\ $h(r)\leqq
\epsilon T_f(r)||_\epsilon$, $h(r)=O(T_f(r))$\hskip1pt).

For a quasi-projective manifold $V$ and for a holomorphic curve
$f:\C \to V$ we write simply $T_f(r)=T_f(r; H)$ for the order function with
respect to an ample line bundle $H$ over a projective compactification
$\bar M$ of $M$ if the choice of $\bar M$ and $H$ do not matter.

The following related property of order functions will be frequently used
(\cite{\NO} Lemma (6.1.5)).

\begin{lemma}
\label{order}
Let $\eta: V \to W$ be a rational mapping between
quasi-projective manifolds $V$ and $W$.
Then for an algebraically non-degenerate holomorphic curve
$f:\C \to V$
$$
T_{\eta\circ f}(r) =O(T_{f}(r)).
$$
Moreover, if $\eta$ is generically finite, then
$$
T_f(r)=O(T_{\eta\circ f}(r)).
$$
\end{lemma}

We define the proximity function $m_f(r; \I)$ not only for
divisors but also for a coherent ideal sheaf $\I$
of the structure sheaf $\O_M$ over $M$.
Let $\{U_j\}$ be a finite open covering of $M$ 
such that
\begin{enumerate}
\item
there is a partition of unity $\{c_j\}$ associated with $\{U_j\}$,
\item
there are finitely many sections
$\sigma_{jk} \in \Gamma(U_j, \I), k=1,2,\ldots$,
generating every fiber $\I_x$ over $x \in U_j$.
\end{enumerate}
Setting $\rho_\I(x) = \left(
\sum_j c_j(x)\sum_k |\sigma_{jk}(x)|^2\right)^{1/2}$,
we take a positive constant $C$ so that
$$
C\rho_\I(x) \leqq 1, \qquad x \in M.
$$
Using the compactness of $M$, one easily verifies that,
up to addition by a bounded continuous function on $M$,
$\log\rho_\I$ is independent of the choices of the open
covering, the partition of unity, the local generators of the
ideal sheaf $\I$, and the constant $C$.

We define the proximity function of $f$ for $\I$ or for
the subspace (may be non-reduced) $Y=(\supp \O_M/\I, \O/\I)$ by
\begin{equation}
\label{2.3}
m_f(r; Y)=m_f(r; \I) =
\int_{|z|=r} \log\frac{1}{C\rho_\I(f(re^{i\theta}))}
\frac{d\theta}{2\pi}\quad (\geqq 0),
\end{equation}
provided that $f(\C) \not\subset \supp Y$.
Note that if $\I$ is the ideal sheaf defined by an effective
divisor $D$ on $M$, $m_f(r; \I)$ coincides $m_f(r; D)$
defined in \cite{NWY02} up to $O(1)$-term.
The function $\rho_\I\circ f(z)$ is smooth over
$\C \setminus f^{-1} (\supp Y)$. 
For $z_0 \in f^{-1} (\supp Y)$
choose an open neighborhood $U$
of $z_0$ and a positive integer $\nu$ such that
$f^*\I=((z-z_0)^\nu)$. Then
$$
\log \rho_\I \circ f(z)= \nu \log |z-z_0|+\psi(z),
\quad z \in U.
$$
for some smooth function $\psi(z)$ defined on $U$.
We define the counting function $N(r; f^*\I)$ and
$N_l(r; f^*\I)$ by using  $\nu$ in the same way as
using $\ord_{z_0}(E)$ in the definition of
$N(r; E)$ and $N_l(r; E)$.
Moreover we define
\begin{align}
\label{current}
\omega_{\I,f}&=\omega_{Y, f}
=- dd^c \psi(z)=- \frac{i}{2\pi} \del\delbar \psi(z)\\
\nonumber
&= dd^c \log \frac{1}{\rho_\I \circ f(z)}
\quad (z \in U),
\end{align}
which is well-defined on $\C$ as a smooth (1,1)-form.
The order function of $f$ for $\I$ or $Y$ is defined by
\begin{equation}
\label{orderideal}
T(r; \omega_{\I,f})=T(r; \omega_{Y, f})
=\int_1^r \frac{dt}{t}\int_{|z|<t} \omega_{\I, f}.
\end{equation}
When $\I$ defines a divisor $D$ on $M$, we see that
$$
T(r; \omega_{\I,f})=T_f(r; L(D))+O(1).
$$

Let $\mathcal{I}_i$ ($i=1,2$) be coherent
ideal sheaves of $\O_M$ and let $Y_i$ be the subspace defined by $\I_i$.
We write $Y_1 \supset Y_2$ if $\I_1 \subset \I_2$.

\begin{theorem}
\label{fmt}
Let $f: \C \to M$ and $\I$ be as above.
Then we have the following:

{\rm (i) (First Main Theorem)}
$$
T(r; \omega_{\I,f})= N(r; f^*\I) + m_f(r; \I) - m_f(1; \I).
$$

{\rm (ii)} If $M$ is projective,
$
m_f(r, \I)=O(T_f(r)).
$

{\rm (iii)}  Let $\mathcal{I}_i$ ($i=1,2$) be coherent
ideal sheaves of $\O_M$ and let $Y_i$ be the subspace defined by $\I_i$.
If $\I_1 \subset \I_2$ or equivalently $Y_1 \supset Y_2$, then
$$
m_f(r; \I_2) \leqq m_f(r; \I_1)+O(1),
$$
or equivalently,
$$
m_f(r; Y_2) \leqq m_f(r; Y_1)+O(1).
$$

{\rm (iv)}  Let $\phi:M_1 \to M_2$ be a holomorphic mappings
between compact complex manifolds.
Let $\I_2\subset \O_{M_2}$  be a coherent ideal sheaf
and let $\I_1\subset \O_{M_1}$ be the coherent ideal sheaf
generated by $\phi^*\I_2$.
Then
$$
m_f(r; \I_1) = m_{\phi\circ f}(r; \I_2) +O(1).
$$

{\rm (v)}  Let $\I_i$, $i=1,2$ be two coherent ideal sheaves of $\O_M$.
Suppose that $f(\C) \not\subset \supp (\O_M/\I_1 \tensor \I_2)$.
Then we have
$$
T(r; \omega_{\I_1\tensor\I_2,f})=T(r; \omega_{\I_1,f})+
T(r; \omega_{\I_2,f})+O(1).
$$
\end{theorem}

\begin{proof}
(i) This immediately follows from
the well-known Jensen formula (cf.\ \cite{\NO} Theorem (5.2.15)).

(ii) Let $Y$ be the subvariety defined by $\I$.
There is an ample divisor $D$ on $M$ such that
$D \supset Y$ (counting multiplicities).
It follows from Theorem (\ref{fmt}) (iii) that
$$
m_f(r; Y) \leqq m_f(r; D)\leqq T_f(r; L(D))=O(T_f(r)).
$$

(iii) (iv) (v) These are immediate by definition.
\end{proof}

\section{General position}

\begin{convention}
\rm Unless explicitly stated otherwise, all varieties, morphisms,
group actions, compactifications, divisors etc. are assumed
to be algebraic.
\end{convention}

\subsection{General position}
Let $A$ be a semi-abelian variety and let $X$ be a
complex algebraic variety (possibly singular) on which $A$ acts:
$$
(a,x) \in A \times X \to a \cdot x \in X.
$$
Let $Y$ be a subvariety embedded into a Zariski
open subset of $X$.

\begin{definition}
{\rm We say that $Y$ is {\it generally positioned in $X$}
if the closure $\bar Y$ of $Y$ in $X$ contains no $A$-orbit.
If the support of a divisor $E$ on a Zariski open subset of
$X$ is generally positioned in $X$, then
$E$ is said to be generally positioned in $X$.}
\end{definition}

Let $\pi:X_1\to X$ be a blow-up of smooth projective manifolds
on which $A$ acts.
Let $D$ be a divisor on $X$ and let $D_1$ be its strict transform.
Then $D_1\sim\pi^*D-E$, where $E$ is an effective divisor with
support contained in the exceptional locus of the blow-up.
If $\pi$ is the blow-up along a smooth connected submanifold $C\subset X$,
then $E$ is empty unless $C\subset D$.

\begin{lemma}
Assume that $D$ is generally positioned in $X$.
Let $\pi:X_1\to X$ be an equivariant blow-up.
Then $D_1=\pi^*D$, i.e., $E$ is empty.
\end{lemma}

\begin{proof}
Since the blow-up is assumed to be equivariant, its center $C$ must
be an invariant subset, i.e., $C$ is a union of $A$-orbits. Now $D$
is assumed to be generally positioned in $X$.
This implies that $D$ contains no $A$-orbit.
Therefore no irreducible component of $C$ is contained in $D$.
\end{proof}

\begin{corollary}
Assume that $D$ is big and generally positioned in $X$.
Then $D_1$ is big, too.
\end{corollary}

\begin{proof}
This is immediate from $D_1=\pi^*D$.
\end{proof}

Unfortunately the assumption of being generally positioned
can not be dropped.
For example, let us consider $X=\P^2(\C)$.
Let $D$ be a line and let $X_1\to X$ be the blow-up of a point $p$
on the line $D$.
Then $X_1$ is a ruled surface.
It admits a fibration $\tau:X_1\to\P^1(\C)$
which arises as follows:  We may identify $\P^1(\C)$ with
$\P(T_p\P^2(\C))$.
Then for $x\in\P^2(\C)\setminus\{p\}$ we set $\tau(x)$ to be the tangent
line at $p$ of the unique line in $\P^2(\C)$ connecting $p$ and $x$.
Now the strict transform $D_1$ of $D$ turns out to be a fiber of $\tau$.
As a fiber of a holomorphic map, it can not be big. However, $D$,
as an effective divisor on $\P^2(\C)$, is big.

To give another example, consider the blow-up of $\P^2(\C)$
in two points $p,q\in D$.
A blow-up decreases the self-intersection number of
a curve by $1$.
Therefore the self-intersection number of the strict
transform $D_2$ of $D$ under this blow-up $X_2\to X$
is a curve with
self-intersection number $-1$. As a consequence we have
$\dim H^0(X_2,L(nD_2))=1$ for all $n\in\N$.

Note that these examples are equivariant for a suitably chosen
action of $A=(\C^*)^2$, but $D$ is not generally positioned
in $\P^2(\C)$.

On the other hand, bigness can only be destroyed, not created via
blow-up. This follows from the following fact:
$D_1=\pi^*D-E$ where $E$ is effective. Thus fixing a section
$\sigma\in H^0(X_1,E)$ we obtain an injection
$$
H^0(X_1,L(nD_1))\stackrel{\alpha}{\hookrightarrow} 
H^0(X_1,L(n\pi^*D))\cong H^0(X,L(nD))\quad ( \forall n\in\N )
$$
given by mapping a section to its tensor product with $\sigma^n$.
Therefore the  Iitaka $D$-dimension can only decrease
(\cite{I71}).

\begin{lemma}
\label{lemma-good}
Let $\pi:X_1\to X$ be an equivariant blow-up,
let $D$ be a divisor on $X$ which is generally positioned
in $X$, and let $D_1$ be its strict transform.
Then $D_1$ is generally positioned in $X_1$, too.
\end{lemma}

\begin{proof}
If $D_1$ would contain an $A$-orbit $\Omega$, we could infer
that $\pi(\Omega)\subset\pi(D_1)=D$. Since $\pi$ is assumed
to be equivariant, this would imply that $D$ contains an $A$-orbit,
namely $\pi(\Omega)$.
\end{proof}

\subsection{Stabilizer}
Let $A$ be a semi-abelian variety such that
\begin{equation}
0 \to T \to A \overset{\pi}{\to} A_0 \to 0,
\end{equation}
where $T \cong (\C^*)^t$ and $A_0$ is an abelian variety.
Let $D$ be a divisor on $A$.
The stabilizer of $D$ is defined by
\begin{equation}
\label{stab}
\StD=\{a\in A:a+D=D\}^0,
\end{equation}
where $\{\cdot\}^0$ denotes the identity component.

\begin{lemma}
\label{a0-ample}
Let $D$ be an effective divisor on $A$ and let $\bar D$ be
its closure in an equivariant compactification $\bar A$ of $A$.
Let $L_0\in\Pic(A_0)$ and let $E$ be an $A$-invariant divisor
on $\bar A$ such that $L(\bar D)\cong L(E)\tensor \pi^*L_0$.
Assume that $\StD$ is contained in $T$.
Then $L_0$ is ample on $A_0$.
\end{lemma}

\begin{proof}
By \cite{NW04} Lemma 5.2 we obtain $c_1(L_0)\geqq 0$.
We may regard $c_1(L_0)$ as a bilinear form on a vector space $V$
which can be interpreted as the Lie algebra $\Lie(A_0)$ or
the dual of cotangent bundle $\Omega^1(A_0)^*$ over $A_0$.
Assume that $L_0$ is not ample.
Then there is a vector $v\in V\setminus\{0\}$
such that $c_1(L_0)|_{\C v}\equiv 0$.
Choose a direct sum decomposition
(orthogonal with respect to $c_1(L_0)$) 
$V=\C v\oplus V'$ and
let $\omega$ be a $(1,1)$-form which is positive on $V'$,
but annihilates $\C v$.
Then $c_1(L)\wedge\omega^{g-1}=0$ where $g=\dim A_0=\dim V$.
Let $\Omega$ be a $(1,1)$-form on $\bar A$ which is
positive along the fibers of $\bar A\to A_0$ as constructed
in \cite{NW03} Lemma 5.1.
Then 
\[
0=
\int_{\bar A}\Omega^s \wedge \pi^*\left(c_1(L_0)\wedge\omega^{g-1}\right)
=\int_D\Omega^s \wedge \pi^*\left(\omega^{g-1}\right)
\]
By construction of $\omega$ this implies that $v$ is everywhere
tangent to $D$. But in this case $v\in\Lie(A_0)$ is in the Lie algebra
of the stabilizer $\StD$. This is a contradiction.
\end{proof}

\begin{proposition}
\label{3.7}
Let $\bar A$ be a smooth equivariant compactification of a semi-abelian
variety $A$.
Let $D$ be an effective divisor on $A$ and let $\bar D$ be
its closure in $\bar A$.
Then the following properties hold.

\begin{enumerate}
\item
$\bar A \setminus A$ is a divisor with only simple normal crossings.
\item
If $\StD=\{0\}$, then $\bar D$ is big on $\bar A$.
\end{enumerate}
\end{proposition}

\begin{proof}
(i)  This is \cite{NW04} Lemma 3.4.

(ii) Due to \cite{NW04} there is a line bundle $L_0$ on $A_0$ and
an $A$-invariant divisor $E$ on $\bar A$ such that $L(\bar D)\cong
L(E)\tensor \pi^*L_0$. By Lemma~\ref{a0-ample} the triviality
of $\StD$ implies the ampleness of $L_0$.

Now consider the $T$-action. Evidently $E$ is $T$-invariant. Since
$T$ acts only along the fibers of $\pi:\bar A\to A_0$, the line bundle
$\pi^*L_0$ is also $T$-invariant. It follows that for every $g\in T$
the pull-back $g^*D$ is linearly equivalent to $D$.\footnote{
Actually $g^*D\sim D$ holds for every $g\in T$ and every
$T \cong (\C^*)^s$-action on a projective manifold.
This can be deduced from the fact that the Picard variety 
of a projective manifold contains
no rational curves.}
Next we define sets $S_x$ for $x\in A$ as follows:
\[
S_x = \cap_{g\in T:g(x)\in D}\,g^*D.
\]
By this definition we know that for every $y\not\in S_x$ there
is a section $\sigma$ in $L(D)$ 
such that $\sigma(x)=0\ne\sigma(y)$. From the
definition it follows furthermore
that $S_x$ is an algebraic subvariety of $A$.
Using the $A$-invariant
trivialization of the tangent bundle $TA \cong  A\times\Lie(A)$
we can identify $T_x(S_x)$ with a vector subspace of $\Lie(A)$.
In this identification we obtain
\[
T_x(S_x) = \cap_{g\in T:g(x)\in D}\,g^*D = \cap_{g\in T:g(x)\in D}\, T_{g(x)}D
= \cap_{y\in\pi^{-1}(\pi(x))\cap D}\, T_y(D).
\]
Thus $T_x(S_x)$ depends only on $\pi(x)$.
Let $F_x=\pi^{-1}(\pi(x))$. Then all the points in $F_x\cap S_x$
have the same tangent space. It follows that $F_x\cap S_x$ is an
orbit under a Lie subgroup of $T$. On the other hand, $F_x\cap S_x$
is an algebraic subvariety. Therefore $F_x\cap S_x$ is an orbit
under an algebraic subgroup of $T$. A priori this subgroup may
depend on the point $x$. However, $T \cong (\C^*)^s$ contains only
countably many algebraic subgroups. For this reason it follows that
this algebraic subgroup must be the same for almost all points $x\in A$.
Thus there is an algebraic subgroup $H\subset T$ such that 
each connected component of $S_x\cap F_x$
is a $H$-orbit for almost all $x\in A$.
But this implies that $D$ is invariant under $H$.
Since $\StD=\{0\}$, $H$ is finite.
Thus $S_x\to A_0$ is generically finite
for almost all $x\in A$. Combined with the ampleness of $L_0$ this
implies that $D$ is big.
\end{proof}

\begin{proposition}
\label{genpos}
Let $Z$ be a reduced subvariety of $A$ and let $\bar Z$ be
its closure in a smooth equivariant compactification $\bar A$ of $A$.
If $\St(Z)=\{0\}$, then there is an equivariant blow-up
$\bar A^\dagger \to \bar A$ such
that the strict transform of $\bar Z$ is generally positioned
in $\bar A^\dagger$.

In particular, there exists a smooth equivariant compactification
of $A$ in which $Z$ is generally positioned.
\end{proposition}

\begin{proof}
We can find an effective reduced divisor $D$ on $A$ such that
$D \supset Z$ and $\StD=\{0\}$.
Thus it suffices to assume that $Z=D$, a divirsor.
Using a result of Vojta (\cite{V99} Theorem 2.4 (2)) we obtain a
(possibly singular) equivariant completion
$\hat{\text\i}:A\hookrightarrow\hat A$ such that
$D$ is generally positioned  in $\hat A$.
Consider the diagonal embedding $j:A\hookrightarrow\bar A\times\hat A$
given by $j=(i,\hat{\text\i})$ and let $\bar A'$ denote the closure
of the image $j(A)$.
Let $\bar A^\dagger \to \bar A'$ be an equivariant
desingularization (cf.\ \cite{Hi64}, \cite{BM97}).
Then the composed map $\bar A^\dagger \to \bar A$ is
a blow-up of $\bar A$.
Considering the natural projection $\bar A^\dagger \to\hat A$,
we conclude, as in Lemma~\ref{lemma-good}, 
that $D$ is generally positioned in $\bar A^\dagger$.
\end{proof}

\begin{proposition}
Let $A$ be a semi-abelian variety, let $A\to\bar A$ be an equivariant
compactification and let $Z$ be a subvariety of $A$.
Then there is an equivariant blow-up $\tilde A \to \bar A$
such that the quotient $\tilde A/\St(Z)$ exists.
\end{proposition}
\begin{proof}
$\St(Z)$ is an algebraic subgroup of $A$. Hence there is
a quotient morphism $q:A\to A/\St(Z)$. 
Let $A/\St(Z) \hookrightarrow Z$ be an $A$-equivariant smooth compactification.
Then $q$ is  a morphism
from an Zariski open subset of $\bar A$ to $Z$ and thus defines
a rational map from $\bar A$ to $Z$.
Now we just blow up $\bar A$ and $Z$ to remove the
indeterminacies and obtain a regular morphism.
Since $q:A\to A/\St(Z)$ is equivariant, it is clear
that the indeterminacies on $\bar A$ are $A$-invariant
subvarieties. Therefore the blow-up can be done
equivariantly.
\end{proof}

\subsection{Finitely many orbits}
We will need the following auxiliary result.
\begin{lemma}
\label{finiteorbit}
Let $A$ be a semi-abelian variety and $A\hookrightarrow\bar A$
a smooth equivariant algebraic compactification.
Then there are only finitely many $A$-orbits in $\bar A$.
\end{lemma}
\begin{proof}
Let $\tau: \C^n \to A$ denote the universal covering.
Then $A=\C^n/\Gamma$, where $\Gamma=\tau^{-1}\{0\}$.
Note that $\Gamma$ generates $\C^n$ as complex vector space.

Let $H$ be an algebraic subgroup of $A$. Then $H$ is a semi-abelian
variety, too. It follows that the connected component $\hat H$
of $\tau^{-1}(H)$
coincides with the complex vector subspace of $\C^n$ generated
by $\hat H\cap\Gamma$. Evidently there are only countably many
finitely generated subgroups of $\Gamma$. It follows that there
are only countably many algebraic subgroups $H$ of $A$.

Let $p$ be a point in $\bar A$ and let $H=A_p$ be its isotropy group.
Let $Ap$ denote the $A$-orbit through $p$.
Let ${\bar A}^H$ denote the fixed point set of $H$-action,
i.e., ${\bar A}^H=\{x\in \bar A : ax=x, \forall a\in H\}$.
Then ${\bar A}^H$ is a closed algebraic subvariety of $\bar A$.
Let $T_p({\bar A}^H)$ be its Zariski tangent space at $p$.
Because $H$ is reductive,
the $H$-action on $T_p(\bar A)$ is almost effective.
On the other hand, because $H$ acts trivially on ${\bar A}^H$,
the action on $T_p({\bar A}^H)$ is likewise trivial.
Therefore there is an almost effective $H$-action
on the quotient vector space $T_p(\bar A)/T_p({\bar A}^H)$.
Since $H$ is abelian, this implies
$\dim H\leqq \dim\left(T_p(\bar A)/T_p({\bar A}^H)\right)$.
 From this we deduce
\[
\dim(Ap)=\dim A-\dim H \geqq \dim X
-\dim\left(T_p(\bar A)/T_p({\bar A}^H)
\right) = \dim T_p({\bar A}^H)
\]
Since $Ap\subset {\bar A}^H$, it follows that ${\bar A}^H$
is smooth at $p$ and $Ap$ is open in ${\bar A}^H$.
In particular, there is an open neighborhood $W$ of $p$ in
$\bar A$ such
that $Ap$ is the only $A$-orbit in $W$ with $H$ as isotropy group.
By virtue of algebraicity it follows that there are only finitely many
$A$-orbits in $\bar A$ with $H$ as isotropy group.

Since there are only countably many algebraic subgroups of $A$,
we obtain as a consequence that there are only countably many
$A$-orbits in $\bar A$.

Thus $A$ is an algebraic group acting on an algebraic variety $\bar A$
with only countably many orbits. This implies that there are
actually only finitely many orbits.
\end{proof}

\subsection{Action}

Let $A$ be a semi-abelian variety and let $\pNc$ be the
complex projective $N$-space.
Then $A$ acts on the product $A \times \pNc$
by the group action of the first factor:
$$
(a, (b,x)) \in A \times (A \times \pNc) \to a \cdot (b,x)=
(a+b, x) \in A \times \pNc.
$$
Let $p : A \times \pNc \to A$ be the first projection.
Let $X$ be an irreducible algebraic subset of $A \times \pNc$
such that $p(X)=A$.
We set
$$
B=\St(X)=\{a \in A; a\cdot X=X\}^0,
$$
and assume that $\dim B >0$.
Set
$C=A/B$.

Taking direct products with $\pNc$, one extends the projection $A\to C$
to $\tau:A\times \pNc\to C \times\pNc$. This is a $B$-principal bundle.
The subvariety $X$ of $A\times\pNc$ is $B$-invariant;
therefore $X=\tau^{-1}(\tau(X))$. It follows that $\tau(X)$ is a closed
subvariety of $C\times\pNc$ which we can regard as the quotient
$X/B$ of $X$ with respect to the $B$-action.
In particular $\pi=\tau|_X: X \to Y=\tau(X)$ is a
$B$-principal bundle such that the
$B$-action on $X$ is simply the principal right action of $B$ for this
bundle structure.

Let $\hat B$ be a smooth equivariant compactification of $B$.
Then we have a relative compactification
$\hat A\to C$ of
$A \to C$ arising as the $\hat B$-bundle associated to the
$B$-principal bundle $A\to C$.
In other words: $\hat A=A\times_B\hat B$ where
$A\times_B\hat B$ denotes the quotient of $A\times\hat B$ with respect
to the equivalence relation for which $(a,b)\sim(a',b')$ if and only if
there exists an element $g\in B$ such that $ag=a'$ and $b=gb'$.
The projection map $p$ extends to $\hat p: \hat A \times \pNc \to \hat A$.
Let $\hat X$ be the closure of $X$ in $\hat A$.
Then $\hat X=X\times_B\hat B$.
The compactness of $\hat B$ implies that the projection map
$\hat\pi:\hat X\to Y$ is proper.

Let $E \subset X$ be an irreducible algebraic subset such that
\begin{equation}
\label{triv}
B \cap \St(E)=\{0\}.
\end{equation}

\begin{proposition}
\label{genpos2}
Let $\hat X$, $X$, $E$, etc.\ be as above. Assume in addition that $E$ is
of codimension one, i.e., a divisor.
Then there is a $B$-equivariant blow-up
$$
\psi: X^\dagger \to \hat X
$$
with center in $\hat X \setminus X$ such that
$X^\dagger$ has a stratification by $B$-invariant strata
$$
X^\dagger=\cup_\lambda \Gamma_\lambda
$$
satisfying the following properties:
\begin{enumerate}
\item
$\Gamma_\lambda \cong X/B_x \quad (x \in \Gamma_\lambda)$
where $B_x=\{b\in B:b \cdot x=x\}$ is the isotropy group at $x$.
\item
The closure of $E$ in $X^\dagger$
contains none of the strata $\Gamma_\lambda$.
\item The open subset $X$ of $X^\dagger$ coincides with one of the
strata $\Gamma_\lambda$.
\end{enumerate}
\end{proposition}

\begin{proof}
Before starting the proof we make a remark:
Since $X\to Y$ is a $B$-principal bundle, we can define
quotient varieties $X/H$ for all algebraic subgroups $H$ of $B$.
Therefore statement (i) of the proposition makes sense.

Now we start the proof.
We will only consider blow-ups $X^\dagger\to\hat X$ which arise in
the following way: We take an equivariant blow-up $B^\dagger\to\hat B$
and define $X^\dagger=X\times_BB^\dagger$. We recall that there are
only finitely many $B$-orbits in $B^\dagger$ (Lemma \ref{finiteorbit})
and that $X\times_BB^\dagger$ is defined as a quotient of
$X\times B^\dagger$.
Let $\{\Omega_\lambda\}_\lambda$ be the family of $B$-orbits in
$B^\dagger$.
Then a stratification $\{\Gamma_\lambda\}_\lambda$ of $X^\dagger$ is
induced as follows: For each $\lambda$ we define $\Gamma_\lambda$
is the image of $X\times\Omega_\lambda$ under the projection
$X\times B^\dagger\to X\times_BB^\dagger=X^\dagger$.
Each of these $B$-orbits $\Omega_\lambda$ can be written as quotient
of $B$ by some closed algebraic subgroup $H_\lambda$:
$$
\Omega_\lambda \cong B/H_\lambda.
$$
Then $H_\lambda$ is the isotropy group of
the $B$-action on $\Gamma_\lambda$
at any point $x\in\Gamma_\lambda$ and $\Gamma_\lambda=X/H_\lambda$.
Thus the stratification
$\{\Gamma_\lambda\}_\lambda$
of $X^\dagger$ has the properties required by (i),
for every choice of an equivariant blow-up $B^\dagger\to\hat B$.

By construction, the open subset $X$ of $X^\dagger$ coincides
with the open $B$-orbit in $B^\dagger$, hence (iii) follows.

Let us now verify that $B^\dagger\to B$ can be chosen in such a way
that property (ii) holds, too.
For $y\in Y$ let $E_y$ be defined as $E_y=\{p\in E:\pi(p)=y\}$.
We observe that $\bar{E_y}=\pi^{-1}(y)\cap\bar E$ for almost all 
$y\in\pi(E)$.
Using \cite{N81}, Lemma 4.1., we infer from \eqref{triv} that for
a generic point $y\in\pi(E)$ the fiber $E_y$ has a discrete
stabiliser with respect to the $B$-action on $X$.
 Thus we may invoke Proposition \ref{genpos} and
deduce that there exists an equivariant blow-up $B^\dagger\to \hat B$
such that $E_y$ is generally positioned in $B^\dagger$.
Let $X^\dagger\to \hat X$ be the associated blow-up of $\hat X$.
Now $E_y$ being generally positioned in $B^\dagger$ implies
that the closure of $E$ in $X^\dagger$ contains none of the strata
$\Gamma_\lambda$.
\end{proof}

\section{Second main theorem for jet lifts}

Let $A$ be a semi-abelian variety of dimension $n$
and let $T$ be the maximal affine
 subgroup of $A$.
Then $T \cong (\C^*)^t$ and there is an exact sequence
of rational homomorphisms
$$
0 \to T \to A \to A_0 \to 0,
$$
where $A_0$ is an abelian variety.
Let $\bar A$ be a smooth equivariant compactification of $A$.
Set $\del A=\bar A\setminus A$ and let $J_k(\bar A, \log \del A)$
be the logarithmic $k$-jet bundle along $\del A$ (cf.\ \cite{N86}).
Then $A$ acts on $J_k(\bar A, \log \del A)$ and there is an
equivariant trivialization
$$
J_k(\bar A, \log \del A)\cong \bar A \times J_{k, A},
$$
where $A$ acts trivially on the second factor $J_{k, A}=\C^{kn}$.
Let $\bar J_{k,A}$ be a projective compactification of $J_{k,A}$.
With the trivial action of $A$ on $\bar J_{k,A}$ 
and the usual action on $A$ (by translations) and $\bar A$
this yields an $A$-equivariant compactification
$$
\bar J_{k}(\bar A, \log \del A)=\bar A \times \bar J_{k,A}
$$
of $J_k(A)$ with an open $A$-invariant subset
$$
\tilde J_k(A)=A \times \bar J_{k,A}.
$$
For example, we may set $\bar J_{k,A}=\P^{nk}(\C)$ or
 $\bar J_{k,A}=(\P^{n}(\C))^k$.
Then $J_k(A)=J_k(\bar A, \log \del A)|_A$ is a Zariski open
subset of $\bar J_{k}(\bar A, \log \del A)$ and
$$
J_k(A) \cong A \times J_{k, A}.
$$
We set
\begin{align*}
J^{\mathrm{reg}}_k(\bar A, \log \del A)&=\left\{
j_k(g) \in J_k(\bar A, \log \del A); j_1(g)\not=0
\right\} \cong \bar A \times J^{\mathrm{reg}}_{k, A},\\
J^{\mathrm{reg}}_k(A) &= J^{\mathrm{reg}}_k(\bar A, \log \del A)|_A
 \cong A \times J^{\mathrm{reg}}_{k, A},
\end{align*}
of which elements are called {\it regular jets}.

Let $f:\C \to A$ be a holomorphic curve and
$J_k(f): \C \to J_k(A)$ be the $k$-jet lift of $f$.
We denote by $X_k(f)$ (resp.\ $\tilde X_k(f)$)
the Zariski closure of the image $J_k(f)(\C)$ in $J_k(A)$
(resp.\ $\tilde J_k(A)$):
\begin{equation}
\label{jetimage}
X_k(f) \subset J_k(A),\qquad
\tilde X_k(f) \subset \tilde J_k(A).
\end{equation}

\begin{theorem} {\rm (Second Main Theorem)}
\label{4.1}\label{smt}
Let $f: \C \to A$ be an algebraically non-degenerate holomorphic curve.
Let $Z$ be a reduced subvariety of $X_k(f)$.
Then there exists a natural number $l_0$
and a compactification $\bar X_k(f)$ of $X_k(f)$
such that for the closure $\bar Z$ of $Z$ in $\bar X_k(f)$
\begin{align}
\label{prox}
m_{J_k(f)}(r; \bar Z) &=S_f(r),\\
\label{wsmt}
T(r; \omega_{\bar Z , J_k(f)}) &\leqq N_{l_0} (r; J_k(f)^* Z)+S_f(r).
\end{align}

In the case of $k=0$ the compactification $\bar A$ of $A$
can be chosen smooth, equivariant, and independent of $f$;
moreover, if $Z$ is a divisor $D$,
 \eqref{prox} and \eqref{wsmt} take the following forms,
respectively:
\begin{align}
\label{prox0}
m_{f}(r; \bar D) &=S_f(r; L(\bar D)),\\
\label{wsmt0}
T_f(r; L(\bar D)) &\leqq N_{l_0}(r; f^*D)+S_f(r; L(\bar D)).
\end{align}
\end{theorem}

\begin{proof}
Since the very basic idea of the proof is the same as that
of the Main Theorem of [NWY03], it will be helpful to confer it.

We extend the subvariety $Z$ to the closure in $\tilde X_k(f)$ which
is denoted by the same $Z$.

We first prove (\ref{prox}) and \eqref{prox0}.
Set $B=\St(X_k(f))$.
Then we have the quotient maps:
\begin{align*} 
q^B &: A \to A/B=C,\\
q^B_k &: J_k(A) \to J_k(A)/B \cong C \times J_{k, A},\\
\tilde q^B_k &: \tilde J_k(A) \to C \times \bar J_{k,A}.
\end{align*}
By \cite{N98} and [NW03] Lemma 2.3
\begin{equation}
\label{logder}
\dim B>0, \quad T_{q^B_k\circ J_k(f)}(r)=S_f(r).
\end{equation}
Setting $\tilde Y_k=\tilde X_k(f)/B$, we have a quotient map:
$$
\tilde\pi_k: \tilde X_k(f) \to \tilde Y_k \subset C \times \bar J_{k,A}.
$$
Let $\bar B$ be a smooth equivariant compactification of $B$.
Define $\hat A$, $\hat X_k(f)$, $\hat Z$, etc.\
as the partial compactifications of
$A$, $\tilde X_k(f)$, $Z$, etc.\ as in subsection 3.4.
We then have proper maps,
\begin{align*}
\hat q^B_k &: \hat A \times \bar J_{k,A} \to C \times \bar J_{k,A},\\
\hat\pi_k =\hat q^B_k|_{\hat X_k(f)} &: \hat X_k(f) \to \tilde Y_k
\subset C \times \bar J_{k,A},
\end{align*}
whose fibers are isomorphic to $\bar B$.

There are two cases, $B \subset \St(Z)$ and
$B \not\subset \St(Z)$, which we consider separately.

(a) Suppose that $B \subset \St(Z)$.
Set $\hat W=\hat\pi_k(\hat Z)=\hat Z/B$.
Then $\hat W$ has at least codimension one in $\tilde Y_k$.
Let $T \cong (\C^*)^t$ be the maximal affine subgroup of $A$
and let $S$ be that of $B$.
Then $S$ is a subgroup of $T$ and there is a splitting,
$T\cong S \times S'$.
Take an equivariant compactification $\bar S'$ of $S'$ and set
$$
\bar A=\hat A \times_{S'} \bar S'.
$$
Then $\bar A$ is an equivariant compactification of $A$ and $\hat A$.
We have an algebraic exact sequence  
$$
0 \to S' \to C \to C_0 \to 0,  
$$
where $C_0$ is an abelian variety, and an equivariant compactification
$\bar C=C \times_{S'} \bar S'$.
Thus $\hat q^B_k$ extends to
$$
\bar q^B_k: \bar A \times J_{k,A} \to \bar C \times J_{k,A},
$$
Let $\bar X_k(f)$ (resp.\ $\bar Y_k$, $\bar W$) be the closure
of $\hat X_k(f)$ (resp.\ $\hat Y_k$, $\hat W$) in
$\bar A \times \bar J_{k,A}$
(resp.\ $\bar C \times \bar J_{k,A}$).
Thus we have the restriction
$$
\bar\pi_k=\bar q^B_k|_{\bar X_k(f)}: \bar X_k (f) \to \bar Y_k.
$$
Note that $\bar\pi _k$ is surjective and
\begin{equation}
\label{jetdiff}
\bar W \not= \bar Y_k.
\end{equation}
It follows from Theorem \ref{fmt} (ii) and (\ref{logder}) that
\begin{align}
\label{e-case}
m_{J_k(f)}(r; \bar Z) &\leqq m_{\bar\pi_k\circ J_k(f)}(r; \bar W)+O(1)\\
\nonumber
&=O(T_{\bar\pi_k \circ J_k(f)}(r))=S_f(r).
\end{align}

(b) Suppose that $B \not\subset \St(Z)$.  We set
$$
B' =B \cap \St(Z),\quad Z'=Z/B',\quad \tilde X_k'(f)=\tilde X_k(f)/B',
\quad A'=A/B',\quad B''=B/B'.
$$
Moreover, we define $W$ as the image of $Z$ under the
quotient $\tilde X_k'(f)\to \tilde X_k'(f)/B''=\tilde Y_k$.
We have the following
commutative diagram and quotient maps:
$$
\begin{array}{ccccc}
Z & \subsetneqq & \tilde X_k(f) & \subset & A \times \bar J_{k,A}\\
 & \scriptstyle{(\codim=1)} & & & \\
\downarrow &  \empty & \downarrow & \empty &
 \downarrow\, \scriptstyle{q^{B'}_k} \\
  \\
Z' & \subsetneqq & \tilde X'_k(f) & \subset & A' \times \bar J_{k,A}\\
 & \scriptstyle{(\codim=1)} & & & \\
\qquad\downarrow\, \scriptstyle{\tilde\pi'_k|_{Z'}} &
 \empty & \quad \downarrow\, \scriptstyle{\tilde\pi'_k}
 & \empty & \downarrow\, \scriptstyle{q^{B''}_k} \\
  \\
W & \subset & \tilde Y_k & \subset & C \times \bar J_{k,A}
\end{array}
$$
Note that
\begin{equation}
\label{sttrivial}
\St(X'_k(f))=B'', \qquad \St(Z') \cap B''=\{0\}.
\end{equation}
Let $\bar B''$ be a smooth equivariant compactification of $B''$.
We have
\begin{align}
\nonumber
\hat A' &=A' \times_{B''} \bar B'',\\
\nonumber
\hat\del A'&=\hat A' \setminus A',\\
\nonumber
\hat X'_k(f) &= \tilde X'_k(f) \times_{B''} \bar B'',\\
\nonumber
\hat Z' &= \bar Z'\quad (\hbox{the closure} \hbox{ of }
Z' \hbox{ in  } \hat X'_k(f)),\\
\nonumber
\hat\del X'_k(f) &=\hat X' _k(f) \setminus \tilde X'_k(f).
\end{align}
Note that the boundary divisor $\hat\del A'$ has only normal crossings
(Proposition \ref{3.7} (i)).
We obtain proper maps
$$
\begin{array}{ccccc}
\hat Z' & \subsetneqq & \hat X'_k(f) & \subset & \hat A' \times \bar J_{k,A}\\
\\
\qquad\downarrow\, \scriptstyle{\hat\pi'_k|_{\hat Z'}} &
 \empty & \quad \downarrow\, \scriptstyle{\hat\pi'_k}
 & \empty & \downarrow\, \scriptstyle{\hat q^{B''}_k} \\
  \\
\hat W & \subset & \tilde Y_k & \subset & C \times \bar J_{k,A}\; ,
\end{array}
$$
where $\hat W=\hat\pi'_k(\hat Z')$.
By Proposition \ref{genpos2} we have a blow-up
$$
\psi: \hat X_k^{\prime \dagger}(f) \to \hat X'_k (f)
$$
with center in $\hat\del X'_k(f) $, the strict transform
$\hat Z^{\prime \dagger}$ of $\hat Z'$ and the boundary
$$
\Gamma=\hat X_k^{\prime \dagger}(f) \setminus \tilde X'_k(f)
$$
with stratification
$\Gamma=\cup_\lambda \Gamma_\lambda$ such that
\begin{align}
\label{quot}
&\Gamma_\lambda \cong \tilde X'_k(f)/\mathrm{Iso}_x(B'')
\quad (x \in \Gamma_\lambda),\\
\label{jetdiff2}
&\Gamma_\lambda \cap \hat Z^{\prime \dagger} \not= \Gamma_\lambda .
\end{align}

Here, if $k=0$, we use Proposition \ref{genpos} in place of
Proposition \ref{genpos2}, and deduce the stated property for $\bar A$.

Let $\psi_{*l}: J_l(\hat X_k^{\prime \dagger}(f), \log \Gamma) \to
J_l(\hat X'_k(f), \log \hat\del X'_k(f))$ be the morphism
naturally induced by $\psi$.
We consider a sequence of morphisms
\begin{align*}
J_l(\hat Z^{\prime \dagger}, \log \Gamma) \subset
& J_l(\hat X_k^{\prime \dagger}(f), \log \Gamma)\:
 \overset{{\psi_{*l}}_{\empty_{\empty}}}{\to}\:
J_l(\hat X'_k(f), \log \hat\del X'_k(f))\\
& \hookrightarrow
J_l(\hat A'\times \bar J_{k,A}, \log( \hat\del A'\times \bar J_{k,A}))\\
& \cong J_l(\hat A', \log \hat\del A') \times J_l(\bar J_{k,A})\\
& \cong \hat A' \times J_l(J_{k, A'}) \times J_l(\bar J_{k,A})\\
& \overset{{\mathrm{proj.}}_{\empty_{\empty}}}{\to}\:
 J_l(J_{k, A'}) \times J_l(\bar J_{k,A}).
\end{align*}
Thus we have a morphism
$$
\beta_l: J_l(\hat X_k^{\prime \dagger}(f), \log \Gamma) 
 \to J_l(J_{k, A'}) \times J_l(\bar J_{k,A}).
$$
Let $p_l: J_l(\hat X_k^{\prime \dagger}(f)) \to \hat X_k^{\prime \dagger}(f)$ be the projection
to the base space.
Henceforth we obtain a proper morphism
$$
\gamma_l= (\hat\pi'_k \circ \psi \circ p_l) \times \beta_l :
J_l(\hat X_k^{\prime \dagger}(f), \log \Gamma) \to \tilde Y_k
\times  J_l(J_{k, A'}) \times J_l(\bar J_{k,A}).
$$
We claim that for some $l_0 \geqq 1$
\begin{claim}
\label{jetdiff3}
\quad\qquad
$\gamma_{l_0}(J_{l_0}(\hat Z')) \not= \gamma_{l_0}(J_{l_0}(\hat X'_k(f))).$
\end{claim}

Assume contrarily that
$\gamma_l(J_l(\hat Z')) = \gamma_l(J_l(\hat X'_k(f)))$
for all $l \geqq 1$.
Then for an arbitrary $z \in \C$
\begin{equation}
\label{contact}
J_l(q^{B'}_1 \circ J_k(f))(z) \in
\gamma_l(J_l(\hat Z^{\prime \dagger}, \log \Gamma)).
\end{equation}
Fix $z_0 \in \C$. Then
$\hat\pi_k \circ J_k(f)(z_0) \in \tilde Y_k$
and we set
$$
\xi_l=J_l(q^{B'}_1 \circ J_k(f))(z_0) \in
\gamma_l(J_l(\hat Z^{\prime \dagger}, \log \Gamma)),
\qquad l \geqq 1.
$$
Set $\Xi_l=\gamma_l^{-1}(\xi_l)$ for $l \geqq 0$.
Then the restriction $p_l|_{\Xi_l}$ is proper and
$p_l|_{\Xi_l}: \Xi_l \to p_l(\Xi_l)$ is an isomorphism.
We set
$$
\Lambda_l=p_l(\Xi_l), \qquad l=1,2, \ldots .
$$
The sequence of $\Lambda_l \supset \Lambda_{l+1}$, $l=1,2,\ldots$
terminates to $\Lambda_\infty=\Lambda_{l_0} = \Lambda_{l_0+1}
=\cdots$ ($\subset \hat X_k^{\prime \dagger}(f)$)
for some $l_0$.
Then $\Lambda_\infty \not=\emptyset$.
If $\Lambda_\infty\cap \tilde{X'_k}(f)\not=\emptyset$, there is an element
$a \in A'$ such that
$$
a\cdot (J_l(q^{B'}_1 \circ J_k(f))(z_0)) \in J_l(Z'),
\quad \forall l\geqq 0.
$$
By the identity principle we deduce that
$a\cdot \tilde X'_k(f) \subset Z'$; this is absurd.

Now assume that $\Lambda_\infty\cap \Gamma \not=\emptyset$.
There is a point $x_0 \in \Lambda_\infty \cap \Gamma$
such that
$$
(x_0, \xi_l) \in J_l(\hat Z^{\prime \dagger})_{x_0}, \qquad l \geqq 1.
$$
Let $\Gamma_{\lambda_0}$ be the boundary stratum containing
$x_0$.
Let $\alpha: \tilde X'_k(f) \to \tilde X'_k(f)/\mathrm{Iso}_{x_0}(B'')
\cong \Gamma_{\lambda_0}$ be the quotient map.
Then there exists an element $a_0 \in A$ such that
$$
a\cdot (\alpha \circ q^{B'}_1 \circ J_k(f)(z)) \in
\Gamma_{\lambda_0}\cap \hat Z^{\prime \dagger}
$$
in a neighborhood of $z_0$ and hence for all $z \in \C$.
Henceforth a contradiction follows from this, (\ref{jetdiff2})
and the image $J_k(f)(\C)$ being Zariski dense in $X_k(f)$.

This proves  Claim \ref{jetdiff3}.

We infer \eqref{wsmt} and \eqref{wsmt0} as in the proof of the
Main Theorem of [NWY02] p.\ 152 (cf.\ [NWY02] (5.12))
with a modification as follows.
Let $\bar X_k(f)=\bigcup_\alpha U_\alpha$ be a finite affine covering,
and let $\sigma_{\alpha \nu}$ be the defining functions of
$Z \cap U_\alpha$.
Then by Claim \ref{jetdiff3} there is a rational function
$\eta$ on $\tilde Y_{k} \times J_{l_0}(J_{k,A'} \times
 J_{l_0}(\bar J_{k,A}))$,
regarded as a rational function on $J_{l_0}(X_k(f))$
such that
\begin{align*}
&\eta \circ J_{l_0}(q_1^{B'}\circ J_k(f))(z) \not\equiv 0,\quad z \in
 \C,\\
&\eta|_{U_\alpha}= \sum_{\nu}\sum_{0\leqq j \leqq l_0}
a_{\alpha \nu j}d^j\sigma_{\alpha \nu},
\end{align*}
where the coefficients $a_{\alpha \nu j}$ are jet diffenrentials
on $U_\alpha$.
Then we have
\begin{align*}
\eta|_{U_\alpha} &= \sum_{\nu} \sigma_{\alpha \nu}
\sum_{0\leqq j \leqq l_0}
a_{\alpha \nu j} \frac{d^j\sigma_{\alpha \nu}}{\sigma_{\alpha \nu}},\\
\left |\eta|_{U_\alpha}\right| &\leqq \left(\sum_\nu |\sigma_{\alpha
 \nu}|^2\right)^{1/2}
\left( \sum_{\nu} \left(
\sum_{0\leqq j \leqq l_0}
|a_{\alpha \nu j}|
\left| \frac{d^j\sigma_{\alpha \nu}}{\sigma_{\alpha \nu}}
\right| \right)^2  \right)^{1/2},\\
\frac{1}{\left(\sum_\nu |\sigma_{\alpha
 \nu}|^2\right)^{1/2}}
&\leqq 
\frac{1}{\left |\eta|_{U_\alpha}\right|}
\left( \sum_{\nu} \left(
\sum_{0\leqq j \leqq l_0}
|a_{\alpha \nu j}|
\left| \frac{d^j\sigma_{\alpha \nu}}{\sigma_{\alpha \nu}}
\right| \right)^2  \right)^{1/2}.
\end{align*}
Therefore we deduce from \eqref{logder} that
\begin{align*}
m_{J_k(f)}(r; \bar Z) &=S_f(r),\\
N(r; J_k(f)^* Z) &= N_{l_0}(r; J_k(f)^* Z)+S_f(r).
\end{align*}
Combining these with the First Main Theorem \ref{fmt},
we obtain \eqref{prox} and \eqref{wsmt}.

Let us now prove the additional statements for the case $k=0$.
In this case we take the quotient, $q: A \to A/\St(Z)$ and we deal with
the holomorphic curve $q \circ f: \C \to A/\St(Z)$ and the subvariety
 $Z/\St(Z)$.  
In this way it is reduced to the case when $\St(Z)=\{0\}$.
Then the compactification of $A$ due to Proposition \ref{genpos}
works for an arbitrary algebraically nondegenerate
$f: \C \to A$.

If $Z$ is a divisor $D$ on A, then Proposition \ref{3.7} (ii)
implies that $D$ is big and we can deduce
\eqref{prox0} with the help of Lemma \ref{2.4a}.
\end{proof}

\section{Higher codimensional subvarieties of $X_k(f)$}

Let $f: \C \to A$ be a holomorphic curve in a semi-abelian variety
$A$.
We use the same notation, $X_k(f)$, $\St(X_k(f))$, etc.\ as
in the previous section.

The purpose of this section is to prove the following.

\begin{theorem}
\label{5.1}
Let $f:\C \to A$ be a holomorphic curve and
let $Z \subset X_k(f)$ be a subvariety of $\codim_{X_k(f)} Z \geqq 2$.
Then there is a compactification $\bar X_k(f)$ such that
for the closure $\bar Z$ of $Z$ in $\bar X_k(f)$
$$
T(r; \omega_{\bar Z J_k(f)}) \leqq \epsilon T_f(r) ||_\epsilon,
\quad \forall \epsilon>0.
$$
In particular,
\begin{equation}
\label{estcodim2}
N(r; J_k(f)^*Z) \leqq \epsilon T_f(r) ||_\epsilon,
\quad \forall \epsilon>0.
\end{equation}
\end{theorem}

{\it Remark.} (i)  For an abelian variety $A$ this was proved by
\cite{Y04}.

(ii)  As a consequence, estimate \eqref{estcodim2} is independent of the
choice of the compactification $\bar X_k(f)$.

It suffices to prove Theorem \ref{5.1} for irreducible $Z$.
Hence, we assume throughout this section
that $Z$ is {\it irreducible}.

Our proof naturally divides into three steps (a)$\sim$(c). 
Before going to discuss the details,
we give an outline of the proof.

(a)  First, we reduce the case to the one that $A$ admits a splitting
$A=B\times C$ where $B$ and $C$ are semi-abelian varieties such that 
\begin{equation}
\label{eqn:int0}
B\subset \St (X_l(f))\quad \text{for all $l \geqq 0$} 
\end{equation}
and the composition of $f$ and the second projection
$q^B :A\to A/B=C$ satisfies
\begin{equation}
\label{eqn:int}
T_{q^B \circ f}(r) =S_{f}(r).
\end{equation}
By this reduction, we may assume that the variety
$X_l(f)$ has splitting $X_l(f)=B\times (X_l(f)/B)$ for all $l \geqq 0$.

We also make a reduction such that the image of $Z$ under
the second projection $\pi_k :X_k(f)\to X_k(f)/B$ has a Zariski
dense image.
Hence by the assumption $\codim_{X_{k}(f)} Z\geqq 2$,
we may assume $\codim_{\pi_k^{-1}(x)}
{Z \cap \pi_k^{-1}(x)}\geqq 2$ for general $x\in X_k(f)/B$.

(b)  The second step is the main part of the proof.
Using the above reduction, we shall construct auxiliary divisors
$F_l\subset \bar{B}\times (X_{k+l}(f)/B)$ for all $l\geqq 0$
with the following properties:
\begin{enumerate}
\item
$(l+1)N_1(r; {J_k(f)}^*{Z}) \leqq N(r; {J_{k+l}(f)}^*{F_l})
+\epsilon T_{f}{(r)} ||_\epsilon, \forall \epsilon >0$:
\item
$T_{J_{k+l}(f)}(r; L({F_l})) \leqq n(l)T_{\gamma \circ f}(r; {D_B})
+\epsilon T_{f}(r; {D}) ||_\epsilon, \forall \epsilon >0$,\par
where $\gamma :A\to B$
is the first projection,
$D$ is an ample line bundle over $\bar A$,
$D_B$ is an ample line bundle over $\bar B$ and
$n(l)$ is a positive integer such that
$\lim_{l\to \infty} n(l)/l= 0$.
\end{enumerate}

(c)  Finally, by (i) and (ii) above we have
\begin{align*}
N_1(r; {J_k(f)}^*{Z}) &\leqq \frac{1}{l+1}N(r; {J_{k+l}(f)}^*{F_l})
+\frac{\epsilon}{l+1} T_{f}(r; {D}) ||_\epsilon \\
&\leqq \frac{n(l)}{l+1}T_{\gamma \circ f}(r; {D_B})
+\frac{\epsilon}{l+1} T_{f}(r; {D}) ||_\epsilon
\end{align*}
for all $\epsilon >0$ and all integer $l\geqq 0$.
Since $n(l)/l\to 0$ ($l\to \infty$), we have
$$
N_1(r; {J_k(f)}^*{Z})\leqq \epsilon (T_{\gamma \circ f}(r; {D_B})
+T_{f}(r; {D})) ||_\epsilon, \quad \forall \epsilon>0.
$$
Since $T_{\gamma \circ f}(r; {D_B}) = O(T_{f}(r; {D}))$,
the proof is completed.

{\bf (a) Reduction.}  Let $f:\C \to A$ be as above.
Let $I_k: \hat X_k(f)\:\: (\hookrightarrow A \times J_{k, A}) \to J_{k, A}$
be the jet projection.
It follows from \cite{N77} (or \cite{NWY02} Lemma 3.8) that
\begin{equation}
\label{lemlog}
T_{I_k \circ J_k(f)}(r)=S_f(r).
\end{equation}

We need the following.

\begin{lem}
\label{log}
Let the notation be as above.
Let $G=\cap _{l \geqq 0} \St({X_l (f)})$ and let
$q^G: A \to A/G$ be the quotient map.
Then
$$
T_{q^G \circ f}(r) = O(T_{I_k \circ J_k(f)}(r))(=S_f(r)).
$$
\end{lem}

\begin{proof}   This is essentially the same as (\ref{logder})
and follows from the jet projection method;
cf.\ \cite{NW03} Lemma 2.4, \cite{NWY02} Lemma 3.8 and
their proofs.
\end{proof}

\begin{lem}
\label{lem:one}
Let $B \subset A$ be a semi-abelian subvariety.
Put $B'=B \cap (\cap _{l \geqq 0} \St({X_l(f)})$.
Let $q^B : A \to A/B$ and $q^{B'}: A \to A/B'$ be quotient mappings.
Then we have
$$
T_{q^{B'} \circ f}(r)= O(T_{q^B \circ f}(r))+S_f(r).
$$
\end{lem}

\begin{proof}
We write $G=\cap _{l \geqq 0} \St({X_l(f)})$.
Taking the natural embedding $A/B' \to (A/B) \times (A/G)$,
we see that
$$
T_{q^{B'}\circ f}(r)=O(T_{q^B \circ f}(r)+ T_{q^G \circ f}(r)).
$$
Thus the claim follows from Lemma \ref{log}.
\end{proof}

\begin{lem}
\label{lem:two}
Let $A$ and $A'$ be semi-abelian varieties with a surjective homomorphism
$p:A \to A'$.
Let $g:\C\to A'$ be a holomorphic curve.
Then we have a holomorphic curve $\hat{g}:\C\to A$
such that $p \circ \hat{g}=g$ and 
$$
T_{\hat{g}}(r) = O(T_{g}(r)).
$$
\end{lem}

\begin{proof}
Set $n=\dim A$ and $n'=\dim A'$.
Let $\varpi: \tilde A\cong \C^{n} \to A$ and $\tilde A'\cong \C^{n'} \to A'$
be the universal covering.
Then there is a surjective linear homomorphism $\tilde p: \tilde A \to \tilde A'$.
Let $\tilde g: \C \to \tilde A'$ be the lifting of $g$.
Let $g(z)=\sum_{j=1}^{n'} g_j(z) e'_j$ with basis $\{e'_j\}$ of $\tilde A'$.
Take a basis $\{e_j\}$ of $\tilde A$ such that $\tilde p (e_j)=e'_j$,
$1 \leqq j \leqq n'$.
Then we set $\hat g(z)=\varpi( \sum_{j=1}^{n'} g_j(z)e_j)$.
It immediately follows from the definition of order functions
(see \cite{NWY02} \S3) that $\hat g$ satisfies the requirement.
\end{proof}

Now we are going to reduce our proof to the case such that
$A=B\times C$ and that $B$ and $C$ are semi-abelian subvarieties
satisfying \eqref{eqn:int0} and \eqref{eqn:int}.
Let $\mathcal{B}$ be the set of all semi-abelian subvarieties
$B\subset A$ such that
$$
T_{q^B \circ f}(r) = S_f(r).
$$
Then since $\cap_{l \geqq 0} \St({X_l(f)}) \in \mathcal{B}$,
we have $\mathcal{B}\not= \emptyset$.
Let $B \in \mathcal{B}$ be a minimal element of $\mathcal{B}$;
i.e., if $B'\subset B$ and $B'\in \mathcal{B}$, then $B'=B$.
If $B_i \in \mathcal{B}, i=1,2$, we deduce from Lemma \ref{lem:one}
that $B_1 \cap B_2 \in \mathcal{B}$.
Thus we get
$$
B \subset \cap_{l \geqq 0} \St({X_l(f)}).
$$
Put $C=A/B$ and let $q^B :A\to C$ be the quotient map.
By Lemma \ref{lem:two} we may take a holomorphic curve
$g:\C \to A$ such that $q^B \circ g=q^B \circ f$ and 
\begin{equation}
\label{eqn:3081}
T_{g}(r) = S_{f}(r).
\end{equation}

We may assume that the Zariski closure of the image $g(\C)$ is
a semi-abelian subvariety $C' \subset A$ (\cite{N77}, \cite{N81}).
Define the semi-abelian variety $\tilde{A}$ by the following pull-back.
\begin{equation*}
\begin{CD}
\tilde{A}@>p_2>> A \\
@Vp_1VV   @VV{q^B}V \\
C' @>{q^B |_{C'}}>>    C
\end{CD}
\end{equation*}
Then $\tilde A=\{(c,a)\in C'\times A:q^B(c)=q^B(a)\}$.
The inclusion map $i:C'\to A$ yields a map
$\tau:C'\to\tilde A$ defined by $\tau(x)=(x,i(x))$.
Note that this morphism $\tau$ 
is a section for $p_1:\tilde{A}\to C'$.
Hence this bundle is trivial, i.e. 
$\tilde{A}\cong B\times C'$ and $\tilde{A}/B=C'$.

Put $\tilde{f}=g\times f:\C \to \tilde{A}$.
Then by \eqref{eqn:3081} we have
\begin{align}
\label{eqn:3093}
T_{f}(r)&=O(T_{\tilde{f}} (r)),\qquad
T_{\tilde f}(r) = O(T_{f}(r)),\\
\label{eqn:3091}
T_{p_1\circ \tilde{f}} (r) &=  S_{\tilde{f}}(r) .
\end{align}
Put 
\begin{equation}
\label{eqn:3095}
B'=B \cap \left( \cap_{l \geqq 0} \St(X_l (\tilde{f}))\right)
\end{equation}
and $p_1':\tilde{A}\to \tilde{A}/B'$ be the quotient map.
By Lemma \ref{lem:one} and \eqref{eqn:3091}, we have
\begin{equation}
\label{eqn:3094}
T_{p_1' \circ \tilde{f}}(r) = S_{\tilde{f}}(r).
\end{equation}
Put $q^{B'}:A\to A/B'$ be the quotient map.
Then we have
\begin{equation}\label{eqn:3092}
T_{q^{B '} \circ f}(r) = O(T_{p_1' \circ \tilde{f}}(r)).
\end{equation}
Hence by \eqref{eqn:3093}, \eqref{eqn:3094} and \eqref{eqn:3092}
we conclude $B'\in \mathcal{B}$.
Since $B$ is minimal in $\mathcal{B}$, we get $B'=B$.
By \eqref{eqn:3095} we have $B \subset \cap_{l\geqq 0}
\St(X_l(\tilde{f}))$.
Let $p_{2,k}: X_k(\tilde f) \to X_k(f)$ be the morphism induced from
$p_2: \tilde A \to A$.  Set
$$
\tilde Z= p_{2,k}^{-1}(Z) \subset X_k(\tilde f).
$$
Note that
$$
N_1(r; J_k(f)^* Z)= N_1(r; J_k(\tilde{f})^* \tilde Z)
$$
and that \eqref{eqn:3093} holds.

For the reduction we need $\codim_{X_k(\tilde f)} \tilde Z \geqq 2$.
By Lemma \ref{lem:one} we see that
$$
B \subset \left(\cap_{l \geqq 0} \St(X_l(f))\right) \cap
\left(\cap_{l \geqq 0} \St(X_l(\tilde f))\right) .
$$
Thus $p_{2,l}: X_l(\tilde f) \to X_l(f)$ is $B$-equivariant, and
induces a morphism
$$
p^B_{2,l}: X_l(\tilde f)/B \to X_l(f)/B.
$$
Let $\pi_l: X_l(f) \to X_l(f)/B$ be the quotient map.
Then it follows from (\ref{eqn:int}) and (\ref{lemlog}) that
\begin{equation}
\label{5.6}
T_{\pi_l \circ J_l(f)}(r)=S_f(r).
\end{equation}

If the image $\pi_k(Z)$ is not Zariski dense in $X_k(f)/B$,
there is a Cartier divisor $H$ on
$X_k(f)/B$ containing $\pi_k(Z)$.
Then, making use of (\ref{5.6}) and the natural embedding
$X_k(f)/B \hookrightarrow (A/B) \times J_{k, A}$ we get
\begin{align}
\label{nondense}
N_1(r; J_k(f)^*Z) &\leqq N(r; (\pi_k\circ J_k(f))^* H)=
O(T_{\pi_k\circ J_k(f)} (r))\\
\nonumber
&= S_f(r).
\end{align}
Therefore the proof of Theorem \ref{5.1} is finished in this case.

We assume that $\pi_k(Z)$ is Zariski dense in $X_k(f)$, and
has a relative dimension at most $\dim B -2$.
Therefore the relative dimension of $\tilde Z \to X_k(\tilde f)/B$
is at most $\dim B -2$, and hence
$\codim _{X_k(\tilde f)} \tilde Z \geqq 2$.

Hence, by replacing
 $A$ by $\tilde{A}$, $C$ by $C'$, $f$ by $\tilde{f}$ and
$Z$ by $p_2^{-1}(Z)$,
we may reduce our problem to the desired situation
\eqref{eqn:int0} and \eqref{eqn:int}.

Therefore we assume the following in the sequel:

\noindent
(i)  Let $B \subset A$
be a semi-abelian subvariety satisfying
\begin{align}
\label{5.2}
B &\subset \cap_{l \geqq 0} \mathrm{St}(X_l (f)),\\
\label{5.3}
T_{q^B\circ f}(r) &=S_f(r),\\
\label{5.4}
A &\cong B \times (A/B),
\end{align}
where $q^B:A \to A/B$ is the quotient map.

\noindent
(ii)  $\pi_k(Z)$ is Zariski dense in $X_k(f)/B$.

{\bf (b)  Auxiliary divisor.}  Let the notation and the assumption be
as above.
Set $C=A/B$.
We have
\begin{equation}
\label{5.5}
A\cong B \times C.
\end{equation}
Then it naturally induces
$$
X_l (f)\cong B\times (X_l (f)/B) \qquad (l \geqq 0).
$$
Let $\bar B$ be an equivariant compactification of $B$ and set
$\hat X_l (f)=\bar B \times (X_l (f)/B)$.
Let
\begin{align*}
\hat\gamma_l &: \hat X_l(f) \to \bar B,\\
\hat\pi_l &: \hat X_l(f) \to X_l(f)/B
\end{align*}
be the natural projections.

We denote by $Z^{\mathrm{ns}}$ the set of non-singular points of $Z$.

\begin{lem}
\label{lem:jyu}
Let $L \to \bar B$ be an ample line bundle.
Then there is a sequence of natural numbers
$n(1),n(2),n(3),\ldots$ satisfying the following:
\begin{enumerate}
\item
$\lim_{l \to \infty} \frac{n(l)}{l}= 0$.
\item
There exist effective Cartier divisors
$F_l\subset \hat{X}_{k+l}(f)$ and line bundles
$M_l$ on $X_{k+l}(f)/B$ such that
$F_l$ is defined by a non-zero element of
$$
H^0(\hat{X}_{k+l}(f),\hat\gamma _{k+l}^*L^{\otimes n(l)}
\otimes (\hat\pi _{k+l})^{*}M_l)
$$
and that for every point $a\in \C$ with
$J_k(f)(a)\in Z^{\mathrm{ns}}$
$$
\ord_{a}{J_{k+l}(f)}^*{F_l}\geqq l+1.
$$
\end{enumerate}
\end{lem}
\begin{proof}
Let $f_B:\C \to B$ be the holomorphic curve defined by
the composition of $f$ and the first projection $A\to B$.
Let $f_C:\C \to C$ be the holomorphic curve defined by
the composition of $f$ and the second projection $A\to C$.  
Then $f_B$ and $f_C$ have Zariski-dense images.
Let $l\geqq 0$ be an integer,
let $p_{k+l,k}: J_{k+l,A} \to J_{k,A}$ be the natural projection,
and let
$$
T\subset J_{k+l}(A)\times C\times J_{k, A} \cong B\times C\times
 J_{k+l, A}\times C\times J_{k, A}
$$
be the Zariski closed subset defined by
$$
T=\{ (b,c,v,c',v')\in B\times C\times J_{k+l,A}\times C\times
 J_{k, A};  b=0,\ c=c',v'=\ p_{k+l,k}(v)\} .
$$
Let $\lambda :B\times C\times J_{k+l,A }\times C\times J_{k, A} \to
 C\times J_{k+l,A}$
be the product of the second projection and the third projection.
We recall the following from \cite{Y04} Proposition 2.1.1.

\begin{lem}
\label{lem:1}\label{lem:112}
There exists a closed subscheme
$\mathcal{T} \subset J_{k+l}(A)\times C\times J_{k, A}$
 with the following properties:
\begin{enumerate}
\item
$\supp \mathcal{T} =T$.
\item
The restriction
$\lambda '=\lambda |_{\mathcal{T}}:\mathcal{T}\to C\times J_{k+l,A}$
is a finite morphism.
Furthermore the restriction of the direct image sheaf
$\lambda'_*(\mathcal{O}_{\mathcal{T}})$ to
$C\times J_{k+l,A}^{\mathrm{reg}}$
is a rank $l+1$ locally free
$\mathcal{O}_{C\times J_{k+l,A}^{\mathrm{reg}}}$-module.
\item
Let $f:\C \to A$ be a holomorphic curve such that $f_B(a)=0$.
Then
$$
\ord_{a}{J_{k+l}(f)}^* {\mathcal{T}_{\rho \circ J_k(f)(a)}}\geqq l+1. 
$$
\end{enumerate}
\end{lem}

Let $r_1:Z^{\dagger}\to \bar{Z}$ be a desingularization
of $\bar{Z}$ such that $r_1$ gives an isomorphism over $Z^\mathrm{ns}$. 
Put $Y_k=X_{k}(f)/B$.
Consider the sequence of morphisms
\begin{equation}
\label{eqn:sqm}
Z^\mathrm{ns} \: \overset{r_0}{\hookrightarrow} \: Z^{\dagger}\: 
\overset{r_{1_{\empty}}}{\to} \: \bar{Z} \: \overset{r_2}{\hookrightarrow} \:
\hat{X}_{k}(f) \: \overset{\hat\pi_{k_{\empty}}}{\to} \: Y_k.
\end{equation}
Here $r_0$, $r_1\circ r_0$ are open immersions and $r_2$ is
a closed immersion.
Put the composition of morphisms to be
$r=\hat\pi_k\circ r_2\circ r_1:Z^{\dagger}\to Y_{k}$.
Let $Y_k^{\mathrm{fl}}$ be a Zariski open subset of $Y_k$
such that $Y_{k}^{\mathrm{fl}}$ is non-singular and the fibers of
$r:Z^{\dagger}\to Y_{k}$ over $Y_k^{\mathrm{fl}}$
are all of the same dimension $\dim Z^{\dagger}-\dim Y_{k}$.
Then the restriction of the family
$r:Z^{\dagger}\to Y_{k}$ to $Y_k^{\mathrm{fl}}$ is a flat family.

Consider the pull back of the sequence of morphisms
\eqref{eqn:sqm} by the natural projection $B\times Y_{k}\to Y_{k}$:
$$
B\times Z^\mathrm{ns} \: \overset{s_0}{\hookrightarrow} \:
B\times Z^{\dagger} \: \overset{s_{1_{\empty}}}{\to} \:
B\times \bar{Z} \: \overset{s_2}{\hookrightarrow} \:
B\times \hat{X}_{k}(f) \: \overset{s_{3_{\empty}}}{\to}
\: B\times Y_{k}.
$$
Again put the composition of these morphisms to be
$s=s_3\circ s_2\circ s_1:B\times Z^{\dagger}\to B\times Y_{k}$.
Then $s$ maps as
$$
s: (a,z) \in B\times Z^{\dagger} \to (a, r(z))\in B\times Y_{k}.
$$
Let $L$ be an ample line bundle on $\bar{B}$ and set
\begin{equation}
\label{phi}
\phi: (a, w) \in B\times \hat{X}_{k}(f) \to a+\gamma _k(w) \in \bar{B}.
\end{equation}
Let $L^{\dagger}_1$ be the line bundle on $B\times Z^{\dagger}$
which is the pull back of $L$ by the composition of morphisms
$$
B\times Z^{\dagger} \: \overset{s_2\circ s_1}{\to} \:
B\times \hat{X}_{k}(f) \: \overset{\phi}{\to} \: \bar{B}.
$$
Since the restriction of $s$ over $B\times Y_k^{\mathrm{fl}}$
(i.e., $s|_{B \times Y_k^{\mathrm{fl}}}:
B\times (Z^{\dagger}|_{Y_k^{\mathrm{fl}}}) \to B\times Y_k^{\mathrm{fl}}$)
is a flat family,
the semi-continuity theorem \cite{H77} p.\ 288 implies
that
there is a Zariski open subset
$U_n\subset B\times Y_k^{\mathrm{fl}}$ $(n>0)$
such that $H^0((B\times Z^{\dagger})|_{P},L^{\dagger \otimes n}_{1,P})$
are all the same dimensional $\C$-vector spaces for $P\in U_n$.
Put this dimension as $G_n$.
Here $(B\times Z^{\dagger})|_{P}$ denotes the fiber of the morphism
$s:B\times Z^{\dagger}\to B\times Y_{k}$ over $P\in B\times Y_{k}$,
and $L^{\dagger \otimes n}_{1,P}$ is the induced line bundle.
Since the intersection ${\cap_{n\geqq 1}} U_n$ is non-empty,
put $(a,w)\in {\cap_{n\geqq 1}} U_n$ and replacing $L$ by
the pull back by the morphism 
$$
B\ni x \mapsto x+a \in B
$$
we may assume $a=0\in B$.

Now for a positive integer $l>0$,
let $\mathcal{T} _l^{\dagger} \subset A \times
J_{k+l,A}\times C\times J_{k,A}$ be the closed subscheme, 
and let $\lambda :\mathcal{T}_l^{\dagger}\to C\times J_{k+l,A}$
be the morphism obtained in Lemma~\ref{lem:1}.
Then $\lambda$ has the following properties;

\begin{enumerate}
\item
$\lambda$ is finite,
\item
the direct image sheaf $\lambda _*\mathcal{O}_{\mathcal{T} _l^{\dagger}}$
is locally generated by $l+1$ elements as
$\mathcal{O}_{C\times J_{k+l,A}}$ module on
$C\times J_{k+l,A}^{\mathrm{reg}}$,
\item
$\lambda$ induces an isomorphism of the underlying topological spaces of
$\mathcal{T}_l^{\dagger}$ and $C\times J_{k+l,A}$.
\end{enumerate}

Since $Y_{k+l}$ is a Zariski closed subset of $C\times J_{k+l,A}$,
we denote $\sigma _{k+l}:Y_{k+l}\to C$ for the composition
with the first projection $C\times J_{k+l,A}\to C$ and
denote $\eta _{k+l}:Y_{k+l}\to J_{k+l,A}$ for the composition
with the second projection.
We have the closed immersion 
\begin{equation}
\label{eqn:503}
B\times Y_{k+l}\times Y_{k}\subset
B\times C\times J_{k+l,A}\times C\times J_{k,A}\cong
A\times J_{k+l,A}\times C\times J_{k,A},
\end{equation}
where the first inclusion is given by
$$
B\times Y_{k+l}\times Y_{k}\ni (b,v,v')\mapsto
(b,\sigma _{k+l}(v),\eta _{k+l}(v),\sigma _{k}(v'),\eta _k(v'))\in
B\times C\times J_{k+l,A}\times C\times J_{k,A}
$$
and the second identification is given by
$$
B\times C\times J_{k+l,A}\times C\times J_{k,A}\ni (b,c,u,c',u')\mapsto
((b,c),u,c',u')\in A\times J_{k+l,A}\times C\times J_{k,A}.
$$
Let $\mathcal{S}_l\subset B\times Y_{k+l}\times Y_{k}$
be the closed subscheme obtained by the pull-back of
$\mathcal{T}_l^{\dagger}$ by \eqref{eqn:503}.
Let $q:\mathcal{S}_l\to Y_{k+l}$ be the composition with
the second projection
$B\times Y_{k+l}\times Y_{k}\to Y_{k+l}$.
We put
$$
Y_{k+l}^{\mathrm{reg}}=Y_{k+l}\cap (C\times J_{k+l,A}^{\mathrm{reg}}),
$$
 which is the Zariski open subset of $Y_{k+l}$.
Then by the above properties of $\lambda$,
we have the corresponding properties for $q$;
\begin{enumerate}
\item
$q$ is finite,
\item
the direct image image sheaf $q _*\mathcal{O}_{\mathcal{S} _l}$
is locally generated by $l+1$ elements as
$\mathcal{O}_{Y_{k+l}}$-module on $Y_{k+l}^{\mathrm{reg}}$,
\item $q$ gives the isomorphism of under lying topological spaces of
$\mathcal{S}_l$ and $Y_{k+l}$.
\end{enumerate}

We consider the following commutative diagram \eqref{eqn:cdm}
obtained by the base change of \eqref{eqn:sqm} with a sequence
of morphisms
$$
\mathcal{S}_l \hookrightarrow  B\times Y_{k+l}\times Y_{k}\to
B\times Y_{k}\to Y_{k}.
$$
Here $B\times Y_{k+l}\times Y_{k}\to  B\times Y_{k}$
is the natural projection:
$$
B\times Y_{k+l}\times Y_{k}\ni (a,w,w')\mapsto (a,w')\in B\times Y_{k}.
$$
\begin{equation}
\label{eqn:cdm}
\begin{CD}
\mathcal{Z}_l^\mathrm{ns} @>>> B\times Y_{k+l}\times Z^\mathrm{ns}@>>>
B\times Z^\mathrm{ns}@>>> Z^\mathrm{ns}\\
@VVu_0V @VVt_0V @VVs_0V @VVr_0V \\
\mathcal{Z}_l^{\dagger} @>>> B\times Y_{k+l}\times Z^{\dagger}@>>>
B\times Z^{\dagger}@>>> Z^{\dagger}\\
@VVu_1V @VVt_1V @VVs_1V @VVr_1V \\
\mathcal{Z}_l @>v'>> B\times Y_{k+l}\times \bar{Z}@>>>
B\times \bar{Z}@>>> \bar{Z}\\
@VVu_2V @VVt_2V @VVs_2V @VVr_2V \\
\cdot @>>> B\times Y_{k+l}\times \hat{X}_{k}(f)@>>>
B\times \hat{X}_{k}(f)@>>> \hat{X}_{k}(f)\\
@VVu_3V @VVt_3V @VVs_3V @VV \hat\pi_k V \\
\mathcal{S} _l @>v>> B\times Y_{k+l}\times Y_{k}@>>>
B\times Y_{k}@>>> Y_{k}
\end{CD}
\end{equation}
Let $\mathcal{L}_l^{\dagger}$ be the line bundle on
$\mathcal{Z} _l^{\dagger}$ obtained by the pull back of
$L^{\dagger}_1$ by the morphisms in the above diagram \eqref{eqn:cdm}.
Let $\mathcal{S}_{l,n}$ be the non-empty Zariski open subset of
$\mathcal{S}_l$ obtained by the inverse image of $U_n$.
Since
$\dim H^0((B\times Z^{\dagger})|_{P},L^{\dagger \otimes n}_{1,P})=G_n$
for $P\in U_n$, the direct image sheaf
$s_{*}L^{\dagger \otimes n}_1$ is a locally free sheaf of
rank $G_n$ on $U_n$ and the natural map
$$
s_{*}L^{\dagger \otimes n}_1\otimes \C  (P)\to
H^0((B\times Z^{\dagger})|_{P},L^{\dagger \otimes n}_{1,P})
$$
is an isomorphism for $P\in U_n$.
This follows by the Theorem of Grauert \cite{H77} p.288,
since $U_n$ is reduced and irreducible.
Here $s:B\times Z^{\dagger}\to B\times Y_{k}$ is the natural map;
i.e., $s=s_3\circ s_2\circ s_1$.
Let $u$ be the morphism $u:\mathcal{Z}_l^{\dagger}\to \mathcal{S}_l$
obtained
by the composition $u=u_3\circ u_2\circ u_1$,
where $u_1$,$u_2$,$u_3$ are the morphisms in the above diagram
\eqref{eqn:cdm}.
Then the natural map 
$$
u_{*}\mathcal{L}_l^{\dagger \otimes n}\otimes \C (P)\to
H^0(\mathcal{Z}^{\dagger}_{l}|_{P},
\mathcal{L}^{\dagger \otimes n}_{l,P})
$$
is also surjective, so an isomorphism on $P\in \mathcal{S}_{l,n}$.
This follows by the Theorem of Cohomology and Base Change
\cite{H77} p.\ 290.
Hence $u_{*}\mathcal{L}_l^{\dagger \otimes n}$ is locally generated by
$G_n$ elements as an $\mathcal{O}_{\mathcal{S}_{l}}$-module on
$\mathcal{S}_{l,n}\subset \mathcal{S}_{l}$.
Let $Y_{k+l,n}=q(\mathcal{S}_{l,n})$ be a non-empty Zariski
open subset of $Y_{k+l}$ (note that the under lying topological
spaces of $\mathcal{S}_{l}$ and $Y_{k+l}$ are the same).
Then by the above properties of $q$, the direct image sheaf
$(q\circ u)_{*}\mathcal{L}_l^{\dagger \otimes n}$
is locally generated by $(l+1)G_n$ elements as a
$\mathcal{O}_{Y_{k+l}}$-module on
$Y_{k+l,n}\cap Y_{k+l}^{\mathrm{reg}}$.
Here, note that $Y_{k+l}^{\mathrm{reg}}$ is non-empty
(otherwise $f$ must be constant) and $Y_{k+l}$ is irreducible.
Hence $Y_{k+l,n}\cap Y_{k+l}^{\mathrm{reg}}$ is also non-empty.

Now look at the following commutative diagram
\begin{equation*}
\begin{CD}
\mathcal{Z}_l^\mathrm{ns}\\
@VVu_0V \\
\mathcal{Z}_l ^{\dagger}@>t_2\circ v'\circ u_1>>
B\times Y_{k+l}\times \hat{X}_{k}(f) @>\psi >>
\bar{B}\times Y_{k+l}@>\rho >> \bar{B}\\
@VVq\circ uV @VV \text{2nd proj} V @VV\tau V \\
Y_{k+l} @= Y_{k+l} @=  Y_{k+l}
\end{CD}
\end{equation*}
where $\rho$ is the first projection,
$\tau$ is the second projection and $\psi$ is the morphism
$$
\psi :B\times Y_{k+l}\times \hat{X}_{k}(f) \ni (a,v,w)
\mapsto (a+\gamma _k(w),v)\in \bar{B}\times Y_{k+l}.
$$
Since
$(\rho \circ \psi \circ t_2\circ v'\circ u_1)^*L=\mathcal{L}_l^{\dagger}$,
we have a natural morphism
\begin{equation}
\label{eqn:530}
\tau _{*}\rho ^{*}L^{\otimes n}=
H^{0}(\bar{B},L^{\otimes n})\otimes _{\C }
\mathcal{O}_{Y_{k+l}} \to (q\circ u)_{*}\mathcal{L}_l^{\dagger \otimes n}.
\end{equation}
Here, note that $\rho \circ \psi =\phi \circ \beta$ where
$\beta :B\times Y_{k+l}\times \hat{X}_{k}(f)\to B\times \hat{X}_{k}(f)$
is the morphism in the diagram \eqref{eqn:cdm} and
$\phi$ was defined by (\ref{phi}).

Now put $I_n=\dim _{\C }H^{0}(\bar{B},L^{\otimes n})$.
Then there is a positive integer $n_0$ and positive constants
$C_1$, $C_2$ such that
$$
I_n>C_1 n^{\dim \bar{B}},\quad G_n<C_2 n^{\dim \bar{B} -2}\quad
\hbox{for } n>n_0.
$$
Here note that
$G_n=\dim _{\C }H^0(B\times Z^{\dagger}|_{P},
L^{\dagger \otimes n}_{1,P})$
for $P\in \cap_{n\geqq 1} U_n$,
and $B\times Z^{\dagger}|_{P}=s^{-1}(P)$ has dimension
$\leqq \dim \bar{B}-2$, for
$\codim_{\hat{X}_{k}(f)} \bar{Z} \geqq 2$ and
$\hat\pi_k\circ r_2:\bar{Z}\to Y_{k}$ is dominant.
Hence for a positive integer $l$,
we can take a positive integer $n(l)$ (e.g.\ $\sim l^{3/4}$) such that
$$
I_{n(l)}>(l+1)G_{n(l)},\qquad
\lim_{l \to \infty} \frac{n(l)}{l}= 0.
$$
Let $\mathcal{F}$ be the kernel of \eqref{eqn:530} for $n=n(l)$;
$$
0\to \mathcal{F} \to \tau _{*}\rho ^{*} L^{\otimes n(l)}\to
(q\circ u)_{*}\mathcal{L}_l^{\dagger \otimes n(l)} \hbox{ (exact)}.
$$
Then we have $\mathcal{F} \not= 0$.
By taking the tensor of a sufficiently ample line bundle $M_l$ on
$Y_{k+l}$
with $\mathcal{F}$, we may assume that
$H^0(Y_{k+l}, \mathcal{F} \otimes M_l)\not= 0$.
Since we have
\begin{equation*}
\begin{split}
H^0(Y_{k+l},\mathcal{F} \otimes M_l)&\subset
H^0(Y_{k+l},(\tau _{*}\rho ^{*}L^{\otimes n(l)})\otimes M_l)\\
&=H^0(Y_{k+l},\tau _{*}(\rho ^{*}L^{\otimes n(l)}\otimes \tau ^*M_l))\\
&=H^0(\bar{B}\times Y_{k+l},\rho ^{*}L^{\otimes n(l)}\otimes \tau ^*M_l),
\end{split}
\end{equation*}
we may take a divisor $F_l\subset \bar{B}\times Y_{k+l}$
which is defined by a non-zero global section of
$H^0(Y_{k+l},\mathcal{F} \otimes M_l)$.
Then we have
$$
\mathcal{Z}_l^\mathrm{ns} \subset  \psi ^{*}F_l.
$$
Here note that $\mathcal{Z}_l^\mathrm{ns} \subset \mathcal{Z}_l$
is an open immersion and
$\mathcal{Z}_l \: \overset{t_2\circ v'}{\hookrightarrow}\:
B\times Y_{k+l}\times \hat{X}_{k}(f)$
is a closed subscheme.

Using the decomposition $A= B\times C$,
we let $f_B:\C \to B$ be the holomorphic curve obtained
by the composition of $f$ and the first projection $A \to B$,
and let $f_C:\C \to C$ be the holomorphic curve obtained
by the composition of $f$ and the second projection $A\to C$.
Now let $a\in \C $ be a point such that $J_k(f)(a)\in Z^\mathrm{ns}$.
Put $\tilde f:\C  \to B\times Y_{k+l}\times \hat{X}_{k}(f)$ as
$$
\tilde f(z)=(f_B(z)-f_B(a),\hat\pi _{k+l} \circ J_{k+l}(f)(z),
J_{k}(f)(a)).
$$
Then we have
$$
\tilde f (\C)\subset B\times Y_{k+l}\times Z,
\quad \tilde f(a)\in \supp\mathcal{Z}_l^\mathrm{ns},
\quad
\psi \circ \tilde f=J_{k+l}(f),
$$
where the last equality holds under the identification
$\bar{B}\times Y_{k+l}= \hat{X}_{k+l}(f)$.\par
Since $v'$ is the base change of $v$ in \eqref{eqn:cdm} and
$\tilde f$ factors through $t_2$, we have
$$
\ord_{a}{\tilde f}^*{\mathcal{Z}_{l}}=
\ord_{a}{(t_3\circ \tilde f)}^*{\mathcal{S}_l},
$$
hence by the construction of $\mathcal{S}_l$ and Lemma \ref{lem:112},
we have
$$
\ord_{a} {\tilde f}^* {\mathcal{Z}_{l}}=
\ord_{a}{(J_{k+l}(f)-f(a))}^*
{\mathcal{T}^{\dagger}_{l,(J_k(f)-f(a))(a)}}\geqq l+1.
$$
Hence we have
$$
\ord_{a}{J_{k+l}(f)}^*{F_l}=
\ord_{a}{\tilde f}^*{\psi ^{*} F_l}\geqq
\ord_{a}{\tilde f}^*{\mathcal{Z}_{l}}^\mathrm{ns}=
\ord_{a}{\tilde f}^*{\mathcal{Z}_{l}}\geqq l+1.
$$
Here note that we consider $F_l$ as the divisor on
$\hat{X}_{k+l}(f)$ by the identification 
$B\times Y_{k+l}\cong \hat{X}_{k+l}(f)$,
and $\tau$ correspond to $\pi _{k+l}$ by this identification.
\end{proof}

{\bf (c) The end of the proof.}
By \eqref{prox} and \eqref{wsmt}  it suffices to show
$$
N_{l_0}(r; J_k(f)^*Z) \leqq 
\epsilon T_{f}(r) ||_{\epsilon}, \quad \forall \epsilon>0.
$$
Note
that
$$
N_{l_0}(r; J_k(f)^*Z) \leqq l_0 N_{1}(r; J_k(f)^*Z).
$$
Furthermore, it suffices to prove
\begin{equation}
\label{eqn:1}
N_1(r; J_k(f)^* Z^{\mathrm{ns}}) \leqq
\epsilon T_{f}(r) ||_{\epsilon}, \quad \forall \epsilon>0.
\end{equation}
For we have
\begin{equation*}
N_1(r; J_k(f)^* Z) =N_1(r; J_k(f)^* Z^{\mathrm{ns}})+
N_1(r; J_k(f)^*(Z \setminus Z^{\mathrm{ns}}))
\end{equation*}
and the second term of the right hand side is estimated to be
at most ``$\epsilon T_f(r)||_\epsilon$''
by induction on dimension of $Z$.
Here note that $\dim Z>\dim (Z\setminus Z^{\mathrm{ns}})$.

It follows from
Lemma \ref{lem:jyu} and (\ref{5.6}) that
\begin{align}
\label{ns}
 (l+1)N_1(r; {J_k(f)}^* {Z^{\mathrm{ns}}}) &\leqq
 N (r; {J_{k+l}(f)}^* {F_l}) \leqq T_{J_{k+l}(f)}(r; L({F_l}))\\
\nonumber
&=n(l)T_{\gamma _{k+l}\circ J_{k+l}(f)}(r; {L})
+T_{\pi_{k+l}\circ J_{k+l}(f)}(r; {M_l})
\\
\nonumber
&\leqq n(l)T_{f_B}(r; {L})+S_{f}(r).
\end{align}
Using $\lim_{l\to \infty} {n(l)}/({l+1)}= 0$ and
$T_{f_B}(r; {L})= O(T_{f}(r; {D}))$,
we obtain \eqref{eqn:1} and our Theorem \ref{5.1}.

\section{Proof of Main Theorem}

(a)  Let the notation be as in the Main Theorem.
The case of $\codim_{X_k(f)}Z \geqq 2$ was finished
by Theorem \ref{5.1}.
Therefore we assume in the rest of this section that
$Z$ is a reduced Weil divisor $D$ on $A$.

Set $B=\St(X_{k+1}(f))$, which has a positive dimension
(cf.\ \eqref{logder}).

\begin{lem}
\label{codim2}
Assume that $D$ is irreducible and $B \not\subset \StD$.
Taking an embedding $X_{k+1}(f) \hookrightarrow J_1(X_{k}(f))$,
we have
$$
\codim_{X_{k+1}(f)}( X_{k+1}(f)\cap J_1(D)) \geqq 2 .
$$
\end{lem}

\begin{proof}
Let $k=0$.
It is first noted that $J_1(A)$ is the holomorphic tangent bundle
$\mathbf{T}(A)$ over $A$, and $X_1(f) \subset \mathbf{T}(A)$.

Assume that $\codim_{X_1(f)}( X_1(f)\cap J_1(D))= 1$.
Let $E$ be an irreducible component of codimension 1 of
$X_1(f)\cap J_1(D)$.
Let $\pi_1: X_1(f) \to A$ be the natural projection.
Then $E$ is an irreducible component of $X_1(f)\cap \pi_1^{-1}(D)$
and $\overline{\pi_1(E)}=D$.

Now $\overline{\pi_1(E)}=D$
combined with $B\not\subset\StD$ implies that $B$
can not stabilize $E$.
Therefore $B \cdot E$ (resp.\ $B \cdot D$) contains an open
subset of $X_1(f)$ (resp.\ $A$). In fact,
since $B$ and $E$ are algebraic,
$B\cdot E$ contains a $B$-invariant 
Zariski open subset $\Omega$ of $X_1(f)$.

Let $p=f(z_0)\in f(\C)$ be a point with the properties:
\begin{enumerate}
\item The orbit
$B\cdot p$ intersects $D\setminus\Sing(D)$ transversely in a point $q$;
\item
$J_1(f)(z_0)\in\Omega$.
\end{enumerate}
Then we choose an analytic 1-dimensional disk $\Delta \subset B$
which contains the unit element $e_B$ of $B$
and we choose a non-empty open subset $U$ of the non-singular part
$D^{\mathrm{ns}}$ of $D$
containing $q$ such that
\begin{enumerate}
\item
the map $\phi:\Delta \times U \hookrightarrow A$ 
given by $\phi(b,u)=b\cdot u$
is an open embedding,
\item
the subbundle
$\cup_{\zeta \in \Delta}\mathbf{T}(\{\zeta\}\times U)
\subset \mathbf{T}(\Delta \times U)$ with
$\mathbf{T}(\{\zeta\}\times U)\cong \mathbf{T}(U)$
gives rise to a holomorphic foliation
on $\phi(\Delta\times U)\subset A$.
\end{enumerate}

Consider $\hat f(z)=b\cdot f(z+z_0)$ with $b\in B$ such that $b\cdot p=q$
and $p=f(z_0)$. Note that $\hat f(0)=b\cdot p=q$
and $J_1(\hat f)(0)=b\cdot J_1(f)(z_0)\in\Omega$.
Then there is an open neighbourhood $W$ of $0$
in $\C$ such that $J_1(\hat f)(z)\in\Omega$
for all $z\in W$. Since $\Omega\subset B\cdot E\subset B\cdot J_1(D)$,
it follows that $\hat f'(z)$ is tangent to the leaves of the
above defined foliation for all $z\in W$.
By the identity principle
this implies $\hat f(\C)=b\cdot f(\C)\subset D$ which is absurd,
since $f$ is algebraically non-degenerate.

The proof for $k \geqq 1$ is similar to the above.
\end{proof}

(b) {\it Proof of the Main Theorem.  }
Let $D=\sum_i D_i$ be the irreducible decomposition.
By making use of Theorem \ref{smt} we have
\begin{align}
\label{6.1}
T(r; \omega_{\bar D , J_k(f)}) &\leqq N_{k_0}(r; J_k(f)^*D)+S_f(r)\\
\nonumber
&\leqq N_1(r; J_k(f)^*D)
 + k_0 \sum_{i<j} N_1(r; J_k(f)^*(D_i \cap D_j))\\
\nonumber
&\quad + k_0 \sum_i N_1(r; J_{k+1}(f)^*J_1(D_i))+S_f(r).
\end{align}
Since $\codim_{X_k(f)} D_i \cap D_j \geqq 2$ for $i\not=j$, it follows
from Theorem \ref{5.1} that
$$
k_0 \sum_{i<j} N_1(r; J_k(f)^*(D_i \cap D_j))\leqq \epsilon
T_f(r)||_\epsilon, \quad
\forall \epsilon >0.
$$
Note that $J_{k+1}(f)^*J_1(D_i)=J_{k+1}(f)^*(X_{k+1}(f)\cap J_1(D_i))$.
If $B \subset \St(D_i)$, then the image of $D_i$ by $X_k(f) \to X_k(f)/B$
is contained in a divisor on $X_k(f)/B$.
Then as in \eqref{e-case} we infer that
$$
N_1(r; J_{k+1}(f)^*J_1(D_i))\leqq N(r; J_k(f)^*D_i) \leqq  S_f(r).
$$

Suppose that $B \not\subset \St(D_i)$.
It follows from Lemma \ref{codim2} and Theorem \ref{5.1} that
$$
N_1(r; J_{k+1}(f)^*J_1(D))\leqq N_1(r; J_{k+1}(f)^*(X_{k+1}(f)\cap J_1(D_i))
\leqq \epsilon T_f(r)||_\epsilon,
\quad \forall \epsilon>0.
$$
Combining these with (\ref{6.1}), we obtain
$$
T(r; \omega_{\bar D , J_k(f)}) \leqq
 N_1(r; f^*D)+\epsilon C T_f(r)||_\epsilon, \quad \forall \epsilon>0,
$$
where $C$ is a positive constant independent of $\epsilon$.
Now the proof of the Main Theorem is completed.  {\it Q.E.D.}

\section{Applications}

(a)  In \cite{G74} M. Green discussed the algebraic degeneracy of a
holomorphic curve $f: \C \to \pnc$ omitting an effective
reduced divisor $D$ on $\pnc$ with normal crossings and of degree
$\geqq n+2$.
He proved the following theorem and conjectured that it would hold
without the condition of finite order for $f$:

\begin{theorem}
\label{green}{\rm (M. Green \cite{G74})}
Let $f:\C \to \P^2(\C)$ be a holomorphic curve of finite order
and let $[x_0,x_1,x_2]$ be the homogeneous coordinate system of
$\P^2(\C)$.
Assume that $f$ omits two lines $\{x_i=0\}, i=1,2$ , and
the conic $\{x_0^2+x_1^2+x_2^2=0\}$.
Then the image $f(\C)$ lies in a line or a conic.
\end{theorem}

Here we answer his conjecture in more general form:

\begin{theorem}
\label{greenconj}
Let $f: \C \to \pnc$ be a holomorphic curve and
let $[x_0,\ldots, x_n]$ be the homogeneous coordinate system
of $\pnc$.
Assume that $f$ omits hyperplanes given by
\begin{equation}
\label{linear}
x_i=0,\qquad 1 \leqq i \leqq n,
\end{equation}
and a hypersurface defined by
$$
x_0^q+\cdots + x_n^q=0, \qquad q \geqq 2.
$$
Then $f$ is algebraically degenerate.
\end{theorem}

\begin{proof}
Let $f(z)=[f_0(z),\ldots, f_n(z)]$ be a reduced representation of
$f$.
Then $f_i(z)$ have no zero for $1 \leqq i \leqq n$.
The assumption implies the existence of an entire function
$h(z)$ such that
$$
f_0^q(z)+\cdots+ f_n^q(z)=e^{h(z)}.
$$
Write the above equation as
$$
\left(f_0(z)e^{-h(z)/q}\right)^q +\cdots+ \left(f_n(z)e^{-h(z)/q}\right)^q=1.
$$
Changing the reduced representation of $f$, we may have that
\begin{equation}
\label{quadr}
f_1^q(z)+\cdots + f_n^q(z) -1= -f_0^q(z).
\end{equation}
Now we take a holomorphic curve into a semi-abelian variety $A=(\C^*)^n$
with the natural coordinate system $(x_1, \ldots, x_n)$ defined by
$$
g: z \in \C \to (f_1(z), \ldots, f_n(z))\in A.
$$
Define a divisor $D$ on $A$ by
$$
x_1^q + \cdots + x_n^q-1=0.
$$
Let $\bar A$ be a equivariant compactification in which $D$ is
generally positioned.
Let $\bar D$ be the closure of $D$ in $\bar A$.
Note that $\StD=\{0\}$ and that $\ord_z g^*D \geqq 2$ for all
$z \in g^{-1}(D)$ by (\ref{quadr}).
Combining this with the Main Theorem ($k=0$), we see that for arbitrary $\epsilon>0$
\begin{align*}
T_g(r; L(\bar D)) &\leqq N_1(r; g^*D)+\epsilon T_g(r; L(\bar D))||_\epsilon \\
&\leqq \frac{1}{q}N(r; g^*D)+ \epsilon T_g(r; L(\bar D))||_\epsilon \\
&\leqq  \frac{1+q\epsilon}{q}T_g(r; L(\bar D))
 ||_\epsilon.
\end{align*}
This leads to a contradiction for $\epsilon <(q-1)/q$.
\end{proof}

{\it Remark.}  The Zariski closure of the image $f(\C)$ can be
more specified in terms of $g$ defined in the above proof.
It follows from [N98] that the Zariski closure of
$g(\C)$ is a translate $X$ of a proper semi-abelian subvariety of $A$
such that $X \cap D=\emptyset$.

(b)  Let $A$ be a semi-abelian variety as above and let
$X \subset J_k(A)$ be an irreducible algebraic subvariety.
We consider the existence problem of an algebraically nondegenerate
entire holomorphic curve $f: \C \to A$ such that $J_k(f)(\C) \subset X$
and $J_k(f)(\C)$ is Zariski dense in $X$.
This is a problem of a system of algebraic differential equations described
by the equations defining the subvariety $X$.

The first necessary condition for the existence of such solution
$f$  is that $\St(X)\not=\{0\}$ (cf.\ \eqref{logder}).
Now we assume the existence of such $f$.
Then we take a big line bundle $L \to X$ and a section
$\sigma \in H^0(X, L)$ which defines a reduced divisor on $X$.
We arbitrarily fix a trivialization
\begin{equation}
\label{trvl}
J_k(f)^*L \cong \C \times \C,
\end{equation}
and regard $J_k(f)^*\sigma$ as an entire function.

\begin{theorem}
\label{deq}
Let the notation be as above.
Then there is no entire function $\psi(z)$ such that every zero of $\psi(z)$
has  degree $\geqq 2$ and
\begin{equation}
\label{deq1}
J_k(f)^*\sigma(z)=\psi(z), \qquad z \in \C.
\end{equation}
In particular, there is no entire function $\psi(z)$ satisfying
\begin{equation}
\label{deq2}
J_k(f)^*\sigma(z)=(\psi(z))^q, \qquad z \in \C,
\end{equation}
where $q \geqq 2$ is an integer.
\end{theorem}

{\it Remark.}  The property given by \eqref{deq1} or \eqref{deq2}
is independent of the choice of the trivialization \eqref{trvl}.

\begin{proof}
Suppose that there is an entire function $\psi(z)$ satisfying
\eqref{deq1} or \eqref{deq2}.
Then it follows that
$$
N_1(r; J_k(f)^*D) \leqq \frac{1}{2} N(r; J_k(f)^*D).
$$
Combining this with the Main Theorem, we infer the following
contradiction:
$$
T_{J_k(f)}(r; L) \leqq \frac{1}{2}T_{J_k(f)}(r; L)+
\epsilon T_{J_k(f)}(r; L) ||_\epsilon.
$$
\end{proof}

\eject

\bigskip
\baselineskip=12pt
\rightline{J. Noguchi}
\rightline{Graduate School of Mathematical Sciences}
\rightline{University of Tokyo}
\rightline{Komaba, Meguro,Tokyo 153-8914}
\rightline{Japan}
\rightline{e-mail: noguchi@ms.u-tokyo.ac.jp}
\bigskip
\rightline{J. Winkelmann}
\rightline{Institut \'Elie Cartan}
\rightline{Universit\'e Nancy I}
\rightline{B.P.239}
\rightline{54506 Vand\oe uvre-les-Nancy Cedex}
\rightline{France}
\rightline{e-mail:jwinkel@member.ams.org}
\bigskip
\rightline{K. Yamanoi}
\rightline{Research Institute for Mathematical Sciences}
\rightline{Kyoto University}
\rightline{Oiwake-cho, Sakyoku, Kyoto 606-8502}
\rightline{Japan}
\rightline{e-mail: ya@kurims.kyoto-u.ac.jp}

\end{document}